\documentclass[12pt]{article}
\usepackage{epsfig}
\usepackage{psfrag}
\usepackage{amssymb}

\newtheorem{thm}{Theorem}[section] 
\newtheorem{lemma}[thm]{Lemma} 
\newtheorem{Def}[thm]{Definition} 
\newtheorem{cor}[thm]{Corollary} 
\newtheorem{prop}[thm]{Proposition} 
\newtheorem{rmk}[thm]{Remark} 

\newtheorem{examp}[thm]{Example}

\newcommand{\ov}{\overline}
\newcommand{\proof}{\noindent {\bf Proof} \hspace{0.2in}} 
\newcommand{\qed}{\hfill\mbox{\raggedright\rule{.07in}{.1in}}
  \vspace{1ex}} 
\newcommand{\dps}{\displaystyle}
\newcommand{\Section}[1]{\section{#1} \setcounter{equation}{0}}

\newcommand{\arraystart}{\renewcommand{\arraystretch}{1.6}}
\newcommand{\arrayfinish}{\renewcommand{\arraystretch}{1.2}}
 
\title{Versal unfoldings for linear retarded functional differential 
equations}
 
\author{Pietro-Luciano Buono\\Centre de Recherches
  Math\'ematiques\\Universit\'e de Montr\'eal\\C.P. 6128, 
Succ. Centre-Ville\\ 
     Montr\'eal, QC H3C 3J7\\CANADA  \and Victor G. LeBlanc\\Department of
  Mathematics and Statistics\\University of Ottawa\\Ottawa, 
ON K1N 6N5\\CANADA}
\date{August 20, 2002}

\begin{document}

\maketitle

\begin{abstract}
We consider parametrized families of linear retarded functional
differential equations (RFDEs) projected onto finite-dimensional invariant
manifolds, and address the question of versality of the resulting
parametrized family of linear ordinary differential equations.
A sufficient criterion for versality is given in terms of readily
computable quantities.  In the case where the unfolding is not versal,
we show how to construct a perturbation of the original
linear RFDE (in terms of delay differential operators)
whose finite-dimensional projection generates a versal unfolding.
We illustrate the theory with several examples, and comment on the
applicability of these results to bifurcation analyses of
{\em nonlinear} RFDEs. 
\end{abstract}

\pagebreak
\Section{Introduction}

Differential equations are used to model a very wide variety of phenomena.  
Frequently, these differential models contain several parameters which are 
often varied or ``tuned'' to describe more accurately the phenomenon under 
study.  Thus, there is considerable interest, from both a pure and an 
applied point of view, to understand how the properties of solutions of 
a parametrized family of differential equations are affected by variation of 
the parameters.  This philosophy is at the core of bifurcation theory.  

One large class of differential equations which are particularly important 
in applications are {\em retarded functional differential equations} (RFDEs) 
\cite{Hal79,HalVL}, which includes the class of ordinary differential 
equations (ODEs), the class of delay differential equations, as well as 
certain types of integro-differential equations, among others.
These equations are used to model various phenomena in fields ranging from 
mathematical biology~\cite{BBL,GM,Kuang,LM} to  
industrial processes~\cite{SA}, and to atmospheric science~\cite{SS}. 
The theory for both linear and nonlinear RFDEs is rather 
well-developped~\cite{Hal79,HalVL}. 
Essentially these equations behave like abstract ODEs on an 
infinite-dimensional (Banach) phase space.  
Thus, many results which are known for ODEs on 
finite-dimensional spaces have analogs in the context of RFDEs.  
For example, in the neighborhood of an equilibrium point of a nonlinear
RFDE, there exists local invariant manifolds (stable, unstable and center 
manifolds) which are tangent to the corresponding invariant subspaces of 
the linearized equations about the equilibrium point, and on which the 
flow near the equilibrium
is either exponentially attracting (stable manifold), exponentially 
repelling (unstable manifold), or non-hyperbolic (center manifold).  
In the context of bifurcation theory, the center manifold reduction of 
the flow is important, since this is where bifurcations (qualitative changes 
in the flow) take place as parameters are varied. In applications, there 
have been many studies (see for example \cite{BC94}) 
of specific RFDEs where stability and 
bifurcations of equilibria were investigated using the center manifold 
reduction theory developped in \cite{Hal79} and \cite{HalVL}.  
Another example of the similarity between ODEs and
RFDEs is in \cite{FMTB,FMH}, where the theory of Poincar\'e-Birkhoff normal 
forms was extended to RFDEs.

One aspect of the bifurcation theory of RFDEs which, surprisingly, has not 
yet been developped is that of extending Arnold's theory \cite{Ar} on versal 
unfoldings of matrices to the case of parameter-dependent linear RFDEs.  
This is the purpose of our present paper. Versal unfoldings of RFDEs 
for certain singularities have been computed for particular classes of 
RFDEs in the study of restrictions on the possible flows on a center manifold
see~\cite{BB,FM96,RLL}. These unfoldings are computed using the normal form
theory of Faria and Magalh$\tilde{\mbox{\rm a}}$es~\cite{FMH}. 
However no attempt is made to give a systematic treatment of unfolding
of linear RFDEs.

We adopt the following strategy. 
We begin with a linear RFDE $\dot{x}(t)={\cal L}_0(x_t)$ (where ${\cal L}_0$ is
a bounded linear functional operator) whose semiflow restricted to a finite-dimensional
subspace is defined by the matrix $B$.  Using a versal unfolding of $B$, we
explicitly construct a parametrized family ${\cal L}(\alpha)$ of bounded linear
functional operators whose finite-dimensional restricted semiflow is defined by a
versal unfolding ${\cal B}(\alpha)$ of the matrix $B$.
In comparison, the realisation of linear ODEs by linear RFDEs 
obtained by Faria and Magalh$\tilde{\mbox{\rm a}}$es~\cite{FMR} provides an
existence result. Using the Hahn-Banach theorem they show the following. 
For any finite dimensional matrix $B$, a necessary and sufficient condition 
for the existence of a bounded linear operator 
${\cal L}_0$ from $C([-\tau,0],{\mathbb R}^{n})$ into
${\mathbb R}^n$
with infinitesimal generator having 
spectrum containing the spectrum of $B$ is that $n$ be larger than or equal to
the largest number of Jordan blocks associated with each eigenvalue 
of $B$. While some of the proofs in our paper have a flavor similar to
the ones in~\cite{FMR}, our results cannot be deduced from realisability
results.

Recall that for a given $c\times c$ matrix $B$ with complex entries, a
{\em $p$-parameter unfolding} of $B$ is a $p$-parameter analytic family 
of matrices ${\cal B}(\alpha)$ such that ${\cal B}(\alpha_0)=B$ for some 
$\alpha_0\in {\mathbb C}^p$.  The unfolding ${\cal B}(\alpha)$ is said to be a 
{\em versal unfolding} of $B$ if for all $q$-parameter unfoldings 
$A(\beta)$ of $B$ (with $A(\beta_0)=B$), there exists an analytic 
mapping $\phi:{\mathbb C}^q\longrightarrow {\mathbb C}^p$ and an analytic 
family of invertible matrices $C(\beta)$ satisfying
\arraystart
\[
\begin{array}{c}
\phi(\beta_0)=\alpha_0\\
C(\beta_0)=I\\
A(\beta)=C(\beta){\cal B}(\phi(\beta))(C(\beta))^{-1}.
\end{array}
\]
\arrayfinish
Thus, a versal unfolding of $B$ is, up to similarity transformations, a 
general analytic perturbation of $B$.  A versal unfolding is said to 
be {\em mini-versal} if the dimension of the parameter space is the smallest 
possible for a versal unfolding.  Of course, the concept of versal unfolding 
of matrices is of importance not only to linear differential equations, but 
to nonlinear differential equations as well, since questions of stability 
and bifurcations of equilibria in nonlinear systems always involve an 
analysis of an associated linearized system.  It is therefore important
to understand the dependence of the associated matrix on system parameters 
(i.e. in the case where the resulting unfolding is not versal, there are 
restrictions on the movement of the eigenvalues, which may influence the 
possible range of dynamics).

In the space $\mbox{\rm Mat}_{c\times c}$
of $c\times c$ matrices with complex entries, let $\Sigma$ denote the
similarity orbit of the matrix $B$.
We will use the following sufficient criterion for versality  
which can be found in~\cite{Ar} and~\cite{Wi}.
\begin{prop}
Let ${\cal B}(\alpha)$ be a $p$-parameter  unfolding of the matrix $B$.  If
\[
\mbox{\rm Mat}_{c\times c}=T\Sigma_{{\cal B}(\alpha_0)}+D_{\alpha}{\cal B}(\alpha_0)\cdot
{\mathbb C}^p,
\]
where $\alpha_0$ is such that ${\cal B}(\alpha_0)=B$ and 
$T\Sigma_{{\cal B}(\alpha_0)}$ denotes the tangent space to $\Sigma$ at
${\cal B}(\alpha_0)$, then ${\cal B}(\alpha)$ is a versal unfolding of $B$.
If the codimension of $T\Sigma_{B(\alpha_0)}$ is equal to $p$, then the 
unfolding is mini-versal.
\label{ArWiprop}
\end{prop}

As motivation for our analysis in this paper, consider a parameter dependent 
nonlinear RFDE 
\begin{equation}
\dot{x}(t)=F(x_t,\alpha)
\label{eq0}
\end{equation}
which has a non-hyperbolic equilibrium, which we assume for the sake of 
simplicity to be at $x=0$, $\alpha=0$.
As previously mentioned, a crucial step in the analysis of the bifurcations
and stability of this equilibrium is to perform a parameter dependent center 
manifold reduction of the RFDE \cite{FMH} in order to obtain a parameter 
dependent nonlinear $c$-dimensional ODE
\begin{equation}
\dot{z}(t)=Bz+G(z,\alpha),
\label{eq1}
\end{equation}
where $B$ is a $c\times c$ constant matrix, and $G$ is {\em nonlinear in 
$z$ and $\alpha$}.  The emphasis is important, since terms which are parameter
dependent and linear in $z$ are contained in the expression for $G$ and they
are the ones which unfold the matrix $B$. Thus, it is natural to investigate 
the versality of this unfolding, for reasons that we have previously 
mentioned. However, the potential difficulty here is that the reduction 
process for RFDEs to finite-dimensional invariant manifolds introduces 
restrictions on the possible types of nonlinearities $G$ which can be 
achieved in (\ref{eq1})\cite{FMR}.  
Thus, there is no {\em a priori} guarantee that a versal unfolding
of the matrix $B$ can be achieved in $G$ by this reduction process, even 
if we have perhaps many parameters (e.g. more than the codimension of the
singular matrix) in the original RFDE. 
It is clear that only terms which are linear in $x_t$ in 
the right-hand side of (\ref{eq0}) 
contribute to terms which are linear in $z$ in the right-hand side of
(\ref{eq1}).  Thus, we 
restrict our attention to the case where (\ref{eq0}) is a parametrized 
family of linear RFDEs.

We have two main results.  The first gives a sufficient condition on the
parametrized linear RFDE (\ref{eq0}) which guarantees that 
the right-hand side of (\ref{eq1}) is
a versal unfolding of the matrix $B$.  The second main result is twofold:
first, we show that despite the previously mentioned
restrictions on $G$ in (\ref{eq1}), it is always possible to realize a
versal unfolding of $B$ by a suitable choice of parameter-dependent
RFDE in (\ref{eq0}); and furthermore, we give a ``canonical'' method of 
computing such a parameter-dependent RFDE in terms of delay differential 
operators. The theory is then illustrated with several examples.

Although our original interest and motivation lies in bifurcation theory
(i.e. $B$ in (\ref{eq1}) has all of its eigenvalues on the imaginary axis),
there is no additional complication in considering a reduction of 
(\ref{eq0}) to a general finite-dimensional invariant manifold (i.e. not
necessarily a center manifold).  Therefore, we develop the theory in this
general context.

As is the case for Arnold's theory of versal unfoldings of matrices, we 
develop our theory by working in complex spaces; since the diagonalization 
theory is much simpler in this context.  Versal unfoldings in the real 
case can be constructed by a decomplexification of the complex unfolding, 
as is done in \cite{Ar} and which we illustrate in Section 7.

\Section{Reduction of linear RFDEs}

Let $C_n=C([-\tau,0],{\mathbb C}^n)$ be the Banach space of continuous 
functions from the interval
$[-\tau,0]$, into ${\mathbb C}^n$ ($\tau>0$) endowed with uniform norm.
We are interested in the linear homogeneous RFDE
\begin{equation}
\dot{z}(t)={\cal L}_0(z_t),
\label{fde}
\end{equation}
where ${\cal L}_0$ is a bounded linear operator from $C_n$ into 
${\mathbb C}^n$.  We write
\[
{\cal L}_0(\varphi)=\int_{-\tau}^{0}\,d\eta(\theta)\varphi(\theta),
\]
where $\eta$ is an $n\times n$ matrix-valued function of bounded variation 
defined on $[-\tau,0]$.  
Let $A_0$ denote the infinitesimal generator of the
semiflow generated by equation (\ref{fde}).
Then it is well-known that the spectrum $\sigma(A_0)$ of $A_0$ 
is equal to the point spectrum of $A_0$,
and $\lambda\in\sigma(A_0)$ if and only if $\lambda$ satisfies the
characteristic equation
\begin{equation}
\mbox{\rm det}\,\Delta(\lambda)=0,\;\;\;\;\;\;\;\;
\mbox{\rm where}\;\;\Delta(\lambda)=\lambda\,I_n-\int_{-\tau}^0\,
d\eta(\theta)e^{\lambda\theta},
\label{char_eq}
\end{equation}
where $I_n$ is the $n\times n$ identity matrix.
We suppose that $\Lambda\subset {\mathbb C}$ is a non-empty finite set 
of eigenvalues of $A_0$, with corresponding generalized $c$-dimensional 
eigenspace $P$. Using adjoint theory, it is known that we can write
\begin{equation}
C_n=P\oplus Q
\label{splitcn}
\end{equation}
where $Q$ is invariant under the semiflow of (\ref{fde}), and invariant
under $A_0$.

Define $C_n^*=C([0,\tau],{\mathbb C}^{n*})$, where ${\mathbb C}^{n*}$ is 
the $n$-dimensional space of row vectors.  We have the adjoint bilinear 
form on $C_n^*\times C_n$:
\begin{equation}
(\psi,\varphi)_n=\psi(0)\varphi(0)-\int_{-\tau}^0\int_0^{\theta}\psi(\xi
-\theta)d\eta(\theta)\varphi(\xi)d\xi.
\end{equation}
We let $\Phi=(\varphi_1,\ldots,\varphi_c)$ be a basis for $P$, and 
$\Psi=\mbox{\rm col}(\psi_1,\ldots,\psi_c)$ be a basis for the dual 
space $P^*$ in $C_n^*$, chosen so that $(\Psi,\Phi)_n$ is the $c\times c$ 
identity matrix, $I_c$.  In this case, we have $Q=\{\varphi\in C_n\,:\,
(\Psi,\varphi)_n=0\}$.  We denote by $B$ the $c\times c$ constant matrix 
such that $\dot{\Phi}=\Phi B$. The spectrum of $B$ coincides with $\Lambda$.
Using the decomposition (\ref{splitcn}), any $z\in C_n$ can be written as
$z=\Phi\,x+y$, where $x\in {\mathbb C}^c$ and $y\in Q$ is
a $C^1$ function.  The dynamics of (\ref{fde}) on $P$ are then given by
\[
\dot{x}=Bx.
\]

\Section{Parametrized families of linear RFDEs}
Consider now a smoothly parametrized family of 
linear RFDEs
of the form
\begin{equation}
\dot{z}(t)={\cal L}(\alpha)(z_t),
\label{fdep}
\end{equation}
where $\alpha\in {\mathbb C}^p$, and ${\cal L}(\alpha_0)={\cal L}_0$ is as in 
(\ref{fde}), for some $\alpha_0\in {\mathbb C}^p$.  
In the sequel, we will assume that a translation has been 
performed in the parameter space ${\mathbb C}^p$ such that $\alpha_0=0$.
Our outline and notation here follows closely that 
of \cite{FMH}. We rewrite (\ref{fdep}) as the system
\begin{equation}
\begin{array}{lll}
\dot{z}(t)&=&{\cal L}_0(z_t)+[{\cal L}(\alpha)-{\cal L}_0](z_t)\\[0.15in]
\dot{\alpha}(t)&=&0.
\end{array}
\label{split_fdep}
\end{equation}
The solutions of this system are of the form $\tilde{z}(t)=(z(t),\alpha(t))^T
\in {\mathbb C}^{n+p}$ (where the superscript $T$ denotes transpose),
the phase space is $\tilde{C}=C_{n+p}=C([-\tau,0],
{\mathbb C}^{n+p})$, and we write (\ref{split_fdep}) as
\begin{equation}
\dot{\tilde{z}}(t)=\tilde{{\cal L}}_0\tilde{z}_t+\tilde{F}(\tilde{z}_t),
\label{augfdep}
\end{equation}
where $\tilde{{\cal L}}_0((u,v)^T)=({\cal L}_0(u),0)^T$ and
$\tilde{F}((u,v)^T)=([{\cal L}(v(0))-{\cal L}_0](u),0)^T$, $u\in C_n$, $v\in C_p$.

Let $\Lambda$, $P$, $Q$, $\Phi$, $\Psi$, $(\,\,,\,\,)_n$ and $B$
be as in the previous section.  Define $\tilde{P}=P\times {\mathbb C}^p$,
$\tilde{Q}=Q\times R$, where $R=\{v\in C_p\,:\,v(0)=0\}$, and consider for 
bases of $\tilde{P}$ and $\tilde{P}^*$, respectively, the columns of the 
matrix $\tilde{\Phi}$ and the rows of the matrix $\tilde{\Psi}$,
\[
\tilde{\Phi}=\left(\begin{array}{cc}
\Phi&0\\
0&I_p\end{array}\right),\;\;\;\;\;\;\;\;\;
\tilde{\Psi}=\left(\begin{array}{cc}
\Psi&0\\
0&I_p\end{array}\right),
\]
which satisfy $(\tilde{\Psi},\tilde{\Phi})_{n+p}=I_{c+p}$.  We have
$\dot{\tilde{\Phi}}=\tilde{\Phi}\tilde{B}$, where $\tilde{B}=\mbox{\rm diag}(B,0)$.
It follows that we have an invariant splitting $\tilde{C}=\tilde{P}\oplus 
\tilde{Q}$.

Let $BC_n$ denote the space of functions from $[-\tau,0]$ to ${\mathbb C}^n$ 
which are uniformly continuous on $[-\tau,0)$ and with a jump discontinuity 
at $0$.  If we define $X_0:[-\tau,0]\longrightarrow {\mathbb C}^n$ by
\[
X_0(\theta)=\left\{\begin{array}{lc}
I_n&\theta=0\\[0.15in]
0&-\tau\leq\theta<0,
\end{array}\right.
\]
then the elements of $BC_n$ can be written as $\xi=\varphi+X_0\mu$, with
$\varphi\in C_n=C([-\tau,0],{\mathbb C}^n)$ and $\mu\in {\mathbb C}^n$, so that
$BC_n$ is identified with $C_n\times {\mathbb C}^n$.  In order to study
(\ref{augfdep}), we need to consider the space $B\tilde{C}=BC_n\times
BC_p$, which can be identified with $\tilde{C}\times {\mathbb C}^{n+p}$.
Define $Y_0:[-\tau,0]\longrightarrow {\mathbb C}^p$ by
\[
Y_0(\theta)=\left\{\begin{array}{lc}
I_p&\theta=0\\[0.15in]
0&-\tau\leq\theta<0,
\end{array}\right.
\]
Let $\pi:BC_n\longrightarrow P$ denote the projection
\[
\pi(\varphi+X_0\mu)=\Phi [(\Psi,\varphi)_n+\Psi(0)\mu],
\]
where $\varphi\in C_n$ and $\mu\in {\mathbb C}^n$.  We consider the
projection $\tilde{\pi}:B\tilde{C}\longrightarrow \tilde{P}$ given by
\[
\tilde{\pi}((\varphi+X_0\mu,\psi+Y_0\nu)^T)=\tilde{\Phi}\left[
\left(\tilde{\Psi},\left[\begin{array}{c}\varphi\\\psi\end{array}\right]
\right)_{n+p}+\tilde{\Psi}(0)\left[\begin{array}{c}\mu\\\nu\end{array}
\right]\right] 
=
\left[\begin{array}{c}
\pi(\varphi+X_0\mu)\\\psi(0)+\nu\end{array}\right].
\]
We now decompose $\tilde{z}$ in (\ref{augfdep}) according to the splitting
\[
B\tilde{C}=\tilde{P}\oplus\mbox{\rm ker}\,\tilde{\pi},
\]
with the property that $\tilde{Q}\subsetneq\,\mbox{\rm ker}\,\tilde{\pi}$, 
and get
\begin{equation}
\begin{array}{rcl}
{\dps\left[\begin{array}{c}\dot{x}\\\dot{\alpha}\end{array}\right]}&=&
{\dps\tilde{B}\left[\begin{array}{c}x\\\alpha\end{array}\right]+
\tilde{\Psi}(0)\tilde{F}\left(\tilde{\Phi}\left[\begin{array}{c}x\\\alpha
\end{array}\right]
+\left[\begin{array}{c}y\\w\end{array}\right]\right)}\\[0.20in]
{\dps\frac{d}{dt}\left[\begin{array}{c}y\\w\end{array}\right]}&=&
{\dps\tilde{A}_{\tilde{Q}^1}\left[\begin{array}{c}y\\w\end{array}\right]+
(I-\tilde{\pi})
\left[\begin{array}{cc}X_0\\Y_0\end{array}\right]
\tilde{F}\left(\tilde{\Phi}\left[\begin{array}{c}
x\\\alpha\end{array}\right]+\left[\begin{array}{c}y\\w\end{array}\right]
\right)},
\end{array}
\label{projfdep}
\end{equation}
where $x\in {\mathbb C}^c$, $\alpha\in {\mathbb C}^p$, 
$y\in Q^1\equiv Q\cap C_n^1$, $w\in R^1\equiv R\cap C_p^1$ 
($C_{n}^1$ and $C_p^1$ denote respectively the subsets of $C_n$ and $C_p$ consisting
of continuously differentiable functions),
and
$\tilde{A}_{\tilde{Q}^1}$ is the operator from 
$\tilde{Q}^1\equiv \tilde{Q}\cap\tilde{C}^1=Q^1\times R^1$ into 
$\mbox{\rm ker}\,\tilde{\pi}$ defined by
\[
\tilde{A}_{\tilde{Q}^1}\left[\begin{array}{c}\varphi\\\psi\end{array}\right]=
\left[\begin{array}{c}\dot{\varphi}\\\dot{\psi}\end{array}\right]+
\left[
\begin{array}{c}
X_0\\Y_0
\end{array}
\right]
\left(
\tilde{{\cal L}}_0\left[\begin{array}{c}\varphi\\\psi\end{array}\right]-
\left[\begin{array}{c}\dot{\varphi}(0)\\\dot{\psi}(0)\end{array}\right]
\right).
\]
If $A_{Q^1}:Q^1\subset\mbox{\rm ker}\,\pi\longrightarrow \mbox{\rm ker}\,
\pi$ is defined by $A_{Q^1}\varphi=\dot{\varphi}+X_0\,[{\cal L}_0\,(\varphi)
-\dot{\varphi}(0)]$, then (\ref{projfdep}) is equivalent to
\[
\begin{array}{rcl}
{\dps\left[\begin{array}{c}
\dot{x}\\\dot{\alpha}\end{array}\right]}&=&
{\dps\left[\begin{array}{c}Bx\\0\end{array}\right]+\left[\begin{array}{c}
\Psi(0)[{\cal L}(\alpha(0)+w(0))-{\cal L}_0](\Phi\,x+y)\\0\end{array}\right]}\\[0.20in]
{\dps\frac{d}{dt}\left[\begin{array}{c}y\\w\end{array}\right]}&=&
{\dps\left[\begin{array}{c}A_{Q^1}y\\\dot{w}-Y_0\dot{w}(0)\end{array}\right]+
\left[\begin{array}{c}(I-\pi)X_0\,[{\cal L}(\alpha(0)+w(0))
-{\cal L}_0](\Phi\,x+y)\\
0\end{array}\right]}.
\end{array}
\]
Since $w\in R$, it follows that $w(0)=0$, so that we get the following 
equations in $BC_n=P\oplus\,\mbox{\rm ker}\,\pi$
\begin{equation}
\begin{array}{rcl}
\dot{x}&=&Bx+\Psi(0)[{\cal L}(\alpha)-{\cal L}_0](\Phi\,x+y)\\[0.20in]
{\dps\frac{d}{dt}y}&=&A_{Q^1}y+(I-\pi)X_0\,[{\cal L}(\alpha)
-{\cal L}_0](\Phi\,x+y),
\end{array}
\label{precmeqs}
\end{equation}
where $x\in {\mathbb C}^c$ and $y\in Q^1$.

\Section{Reduction to parameter dependent invariant manifold}
In this section, we show that (\ref{projfdep}) admits a local,
semiflow invariant, $c+p$-dimensional manifold in
$B\tilde{C}$, which is tangent at the origin to $\tilde{P}$, 
and such that the dynamics of (\ref{projfdep}) restricted to this 
manifold are linear in $x\in {\mathbb C}^c$.

We want the nontrivial part of our invariant manifold to be of the form 
\begin{equation}\label{cmexpr}
u=\Phi x+h(\alpha)x
\end{equation}
where $h:{\mathbb C}^{p}\rightarrow \mbox{\rm Mat}_{1\times c}(Q^1)$ is 
a smooth map, and $\mbox{\rm Mat}_{1\times c}(Q^1)$ denotes the space of 
$1\times c$ matrices whose elements are in the space $Q^1$ which has been
defined in the previous section.
Now, as an infinite-dimensional system, the RFDE~(\ref{fdep}) can be
written as~\cite{FMTB}
\begin{equation}\label{cmeq1}
\displaystyle\frac{du}{dt}=A u+X_0 [{\cal L}(\alpha)-{\cal L}_0]u.
\end{equation}
Combining~(\ref{cmexpr})
and the first equation in (\ref{precmeqs}) we obtain
\begin{equation}\label{cmeq2}
\begin{array}{rcl}
\displaystyle\frac{du}{dt}&=&(\Phi+h(\alpha))\dot x\\
&=&(\Phi+h(\alpha))(Bx+\Psi(0)[{\cal L}(\alpha)-{\cal L}_0](\Phi+h(\alpha))x).
\end{array}
\end{equation}
The expression~(\ref{cmexpr}) represents a locally semiflow invariant
manifold of~(\ref{split_fdep}) near the origin $(z,\alpha)=(0,0)$ if 
equation~(\ref{cmeq1}) is equal to~(\ref{cmeq2}) for values 
of $\alpha$ in a small neighborhood of $0$. Rearranging the two 
equations and simplifying implies that we need to solve 
\begin{equation}\label{cmeq}
\begin{array}{l}
(A_{Q^{1}}+(I-\pi)X_0 [{\cal L}(\alpha)-{\cal L}_0])h(\alpha)
+(I-\pi)X_0 [{\cal L}(\alpha)-{\cal L}_0](\Phi)\\[0.2in]
-h(\alpha)(B+\Psi(0)[{\cal L}(\alpha)-{\cal L}_0](\Phi+h(\alpha)))=0
\end{array}
\end{equation}
for small values of $\alpha$ in a neighborhood of $0$.

Our main tool for proving this result is the implicit function 
theorem (IFT). Because of the smoothness properties required for the 
application of the IFT, we need to work with the $C^1$ norm instead of 
the uniform norm on the space $Q^1$. We claim that the spectral properties 
of the operators $A_0$, $A$ and $A_{Q^1}$ which are given in 
lemmas 5.1 and 5.2 of \cite{FMTB} remain valid if we replace the uniform 
norm by the $C^1$ norm in the domain of these operators.  In particular, we
have
\begin{equation}
\sigma(A_{Q^1})=P\sigma(A_{Q^1})=\sigma(A_0)\setminus\Lambda,
\label{specpropAs}
\end{equation}
where $\sigma$ denotes the spectrum, and $P\sigma$ is the point spectrum.
Clearly the point spectrum does not depend on the norm.  Also, the resolvent 
set does not change, since the above operators are bounded when we replace 
the uniform norm by the $C^1$ norm in the domain of the operators, and
for every $\lambda$ in the resolvent sets (with respect to the uniform 
norm) of these operators, the resolvent operator 
$({\cal A}-\lambda\,I)^{-1}$ (where ${\cal A}$ is any of these
operators) is defined 
on all of $BC_n$, and not just a dense proper subset.  Therefore,
with respect to the $C^1$ norm, ${\cal A}-\lambda\,I$ is bounded and
surjective, so $({\cal A}-\lambda\,I)^{-1}$ is bounded, by Banach's theorem.
 
\begin{prop}
Let $(Q^1)^c=Q^1\times\cdots\times Q^1$ ($c$ times) endowed with $C^1$ 
norm
\[
||(h^1,\ldots,h^c)||_{(Q^1)^c}=
\sum_{i=1}^c\,(|h^i|_C+|\dot{h^i}|_C),
\] (where $|\;\;|_C$ denotes uniform
norm), and let
$(\mbox{\rm ker}\,\pi)^c=\mbox{\rm ker}\,\pi\,\times\cdots\times
\mbox{\rm ker}\,\pi$ ($c$ times) endowed with norm
\[
||(\psi^1+X_0\,\alpha^1,\ldots,\psi^c+X_0\,\alpha^c)||_{\mbox{\rm 
(ker\,$\pi$)$^c$}}=\sum_{i=1}^c\,(|\psi^i|_C+|\alpha^i|_{{\mathbb C}^n}).
\]
Consider the nonlinear operator $N:{\mathbb C}^p\times (Q^1)^c\longrightarrow
(\mbox{\rm ker}\,\pi)^c$
defined by
\begin{equation}
\begin{array}{lll}
N(\alpha,h)&=&\left(A_{Q^1}+(I-\pi)X_0\,[{\cal L}(\alpha)-{\cal L}_0]\right)h+
(I-\pi)X_0\,[{\cal L}(\alpha)-{\cal L}_0](\Phi)-\\[0.15in]&&
h\,\left(B+\Psi(0)\,[{\cal L}(\alpha)-{\cal L}_0](\Phi+h)\right),
\end{array}
\label{homeq}
\end{equation}
where the actions of $A_{Q^1}$ and ${\cal L}(\alpha)$ on $h\in (Q^1)^c$ are
defined componentwise in the obvious way.
Then, 
there exists a neighborhood $V$ of $0$ in ${\mathbb C}^p$ and a unique
smooth mapping $\alpha\longmapsto h(\alpha)$ from $V$ into $(Q^1)^c$ such that
$h(0)=0$ and such that $N(\alpha,h(\alpha))=0$ for all $\alpha\in V$.
\end{prop}
\proof
It is clear that with the chosen topologies, $N$ is a smooth mapping.  
Moreover, it is also clear that $N(0,0)=0$. The partial Fr\'echet 
derivative $N_h(0,0)$ is the bounded linear operator from $(Q^1)^c$ 
into $(\mbox{\rm ker}\,\pi)^c$ defined by
\[
N_h(0,0)v\equiv Jv=A_{Q^1}v-vB.
\]
We suppose that we have chosen a system of coordinates such that the 
$c\times c$ matrix $B$ is in Jordan canonical form
\[
B=\left(\begin{array}{cccccc}
\lambda_1&\sigma_1&&&&\\
&\lambda_2&\sigma_2&&&\\
&&\cdot .&\cdot .&&\\
&&&\cdot .&\cdot .&\\
&&&&\cdot .&\sigma_{c-1}\\
&&&&&\lambda_{c}
\end{array}\right),
\]
where $\sigma_i=1$ or $\sigma_i=0$.
Suppose that $v=(v^1,\ldots,v^c)\in (Q^1)^c\setminus \{0\}$ is such that
$Jv=0$.  Then this is equivalent to 
\[
A_{Q^1}(v^1,\ldots,v^c)=(\lambda_1\,v^1,\sigma_1v^1+\lambda_2v^2,\ldots,
\sigma_{c-1}v^{c-1}+\lambda_c v^c),
\]
which implies that one of the $\lambda_i$ must be in the point
spectrum of $A_{Q^1}$, which is a contradiction (see (\ref{specpropAs})).  
Therefore, $\mbox{\rm ker}\,J=\{0\}$.  

We now show that $J$ is surjective.  We use an approach similar to the 
proof of Theorem 5.4 of \cite{FMTB}. Let $i\in \{1,\ldots,c\}$, 
$\varphi\in Q^1$, and $\xi\in\mbox{\rm ker}\,\pi$ be such that 
$\xi=(A_{Q^1}-\lambda_i\,I)\varphi$.  Define $v=(v^1,\ldots,v^c)\in (Q^1)^c$
by $v^k=0$, $k\neq i$, and $v^i=\varphi$.  Then $f=Jv$ has the form
$f=(f^1,\ldots,f^c)$, where
\begin{equation}
f^k=\left\{\begin{array}{ll}
\xi&\mbox{\rm if $k=i$}\\
-\sigma_i\varphi&\mbox{\rm if $k=i+1$ and $i<c$}\\
0&\mbox{\rm otherwise}.
\end{array}\right.
\label{fdef}
\end{equation}
Let $g=(g^1,\ldots,g^c)\in (\mbox{\rm ker}\,\pi)^c$.  We now use an
induction argument to show that there exists an 
$h=(h^1,\ldots,h^c)\in (Q^1)^c$ which is such that $Jh=g$.
We know from (\ref{specpropAs}) that $\lambda_1$ is in
the resolvent set $\rho(A_{Q^1})$.  Therefore,
there exists $\varphi_1\in Q^1$ such that $(A_{Q^1}-\lambda_1\,I)\varphi_1
=g^1$. If we define ${\cal H}_1=(\varphi_1,0,\ldots,0)\in (Q^1)^c$, 
then we get from (\ref{fdef}) that $(J{\cal H}_1)^1=g^1$.  Suppose now 
that $i\in\{1,\ldots,c-1\}$ is such that there exists ${\cal H}_i\in (Q^1)^c$
satisfying
\begin{equation}
(J{\cal H}_i)^p=g^p,\;\;\;\;\;\forall p\leq i.
\label{indarg}
\end{equation}
Define $H^{i+1}=(J{\cal H}_i)^{i+1}\in\mbox{\rm ker}\,\pi$.  Since
$\lambda_{i+1}\in\rho(A_{Q^1})$, there exists $\varphi_{i+1}\in Q^1$ such
that $(A_{Q^1}-\lambda_{i+1}\,I)\varphi_{i+1}=g^{i+1}-H^{i+1}$.
If we define ${\cal G}_{i+1}=(\zeta^1,\ldots,\zeta^c)\in (Q^1)^c$ by
$\zeta^{k}=0$ if $k\neq i+1$, $\zeta^{i+1}=\varphi_{i+1}$, it follows
from (\ref{fdef}) that
\[
(J{\cal G}_{i+1})^k=\left\{
\begin{array}{ll}
g^{i+1}-H^{i+1}&\mbox{\rm if $k=i+1$}\\
-\sigma_i\varphi_{i+1}&\mbox{\rm if $k=i+2$ and $i+1<c$}\\
0&\mbox{\rm otherwise}.
\end{array}
\right.
\]
If we now set ${\cal H}_{i+1}={\cal H}_i+{\cal G}_{i+1}$, we get from 
(\ref{indarg})
\[
(J{\cal H}_{i+1})^p=\left\{
\begin{array}{ll}
g^{i+1}&\mbox{\rm if $p=i+1$}\\
g^p&\mbox{\rm if $p\leq i$}.
\end{array}
\right.
\]
By induction on $i$, we have that (\ref{indarg}) holds for $i=c$, and thus 
$J$ is a surjection. It now follows that $J$ has a bounded inverse, and 
we get the conclusion of the Proposition by virtue of the implicit function 
thereom.
\hfill\qed

\begin{prop}
Let $h(\alpha)$ be the solution of $N=0$ in~(\ref{homeq}) defined for all
$\alpha\in V$, where $V$ is some neighborhood of the origin in
${\mathbb C}^p$.  Then the following set
\begin{equation}
W=\{((x,\alpha),(y,w))\in {\mathbb C}^{c+p}\times \tilde{Q}^1 
\,:\,w=0,\,y=h(\alpha)x,\;\;x\in {\mathbb C}^c, \alpha\in V\}
\label{cmdef}
\end{equation}
is a locally semiflow invariant manifold for the system (\ref{projfdep}).  
The (non-trivial) dynamics on this manifold reduce to
\begin{equation}
\dot{x}=Bx+\Psi(0)[{\cal L}(\alpha)-{\cal L}_0](\Phi+h(\alpha))x.
\label{cmeqsparam}
\end{equation}
\end{prop}
\proof
It is obvious that $W$ is tangent to $\tilde{P}$.  The semiflow 
invariance of $W$ follows from substitution of $y=h(\alpha)x$
into (\ref{precmeqs}), and using the fact that $h(\alpha)$ is a solution to
$N=0$ in (\ref{homeq}).
\hfill\qed

\Section{Sufficient condition for versality}
From the previous section, we have seen that the dynamics of
(\ref{fdep}) near the equilibrium solution $(z,\alpha)=(0,0)$ 
reduces to the $c$-dimensional parametrized linear system 
\begin{equation}
\dot{x}={\cal B}(\alpha)x, 
\end{equation}
where by (\ref{cmeqsparam}) we have 
\begin{equation}
{\cal B}(\alpha)=B+\Psi(0)[{\cal L}(\alpha)-{\cal L}_0](\Phi+
h(\alpha)).
\label{Bunfdef}
\end{equation}
Note that ${\cal B}(0)=B$ and since $h(0)=0$, we have
\[
D_{\alpha}{\cal B}(0)=\Psi(0)\left.D_{\alpha}[{\cal L}(\alpha)
(\Phi)]\right|_{\alpha=0}.
\]

Let $\mbox{\rm Mat}_{c\times c}$ denote the space of $c\times c$ matrices 
with complex entries. 
For $i,j\in\,\{1,\ldots,c\}$,
let $E_{ij}\in \mbox{\rm Mat}_{c\times c}$ be the 
matrix whose elements are all 0 except the element in row $i$ and column 
$j$, whose value is 1.  For $U_1$, $U_2\in 
\mbox{\rm Mat}_{c\times c}$, denote $[U_1,U_2]\equiv U_1U_2-U_2U_1$.

Denote the linear mapping $\Theta:\mbox{\rm Mat}_{c\times c}\longrightarrow 
{\mathbb C}^{c^2}$ by $\Theta(E_{ij})=e_{(i-1)c+j}$, where $e_{\ell}$ 
denotes the row vector whose components are all 0 except the 
$\ell^{th}$ which is 1.
\begin{thm}
The parametrized family (\ref{fdep}) generates a versal 
unfolding ${\cal B}(\alpha)$ (see (\ref{Bunfdef}))
of the matrix $B=(\Psi,\dot{\Phi})_n$ if the $(c^2+p)\times c^2$
matrix
\[
S=\left(\begin{array}{ccc}
\hspace*{0.1in}&\Theta([B,E_{11}])&\hspace*{0.1in}\\
\hspace*{0.1in}&\Theta([B,E_{12}])&\hspace*{0.1in}\\
\hspace*{0.1in}&\vdots&\hspace*{0.1in}\\
\hspace*{0.1in}&\Theta([B,E_{cc}])&\hspace*{0.1in}\\[0.15in]
\hspace*{0.2in}&\Theta\left(\left.{\displaystyle 
\Psi(0)\,\frac{\partial}{\partial\alpha_1}\,
[{\cal L}(\alpha)(\Phi)]\,}\right|_{\alpha =
0}\right)&\hspace*{0.2in}\\
\hspace*{0.2in}&\vdots&\hspace*{0.2in}\\
\hspace*{0.2in}&\Theta\left(\left.{\displaystyle 
\Psi(0)\,\frac{\partial}{\partial\alpha_p}\,
[{\cal L}(\alpha)(\Phi)]\,}\right|_{\alpha =
0}\right)&\hspace*{0.2in}
\end{array}\right)
\]
has rank $c^2$.  If, in addition, 
\begin{equation}
\mbox{\rm dim}\,\,\mbox{\rm span}(\Theta([B,E_{11}]),
\Theta([B,E_{12}]),\ldots,\Theta([B,E_{cc}]))=c^2-p,
\label{minitransv}
\end{equation}
then the versal unfolding ${\cal B}(\alpha)$ is mini-versal.
\label{thmsuff}
\end{thm}
\proof
In $\mbox{\rm Mat}_{c\times c}$, the tangent space of the similarity orbit
through $B$ is given by \cite{Ar}
\[
{\cal O}=\{\,[B,X]\,:\,X\in\mbox{\rm Mat}_{c\times c}\,\}.
\]
Thus, if the matrix $S$ has rank $c^2$, we have
\[
\mbox{\rm Mat}_{c\times c}={\cal O}+
\Psi(0)\left.D_{\alpha}[{\cal L}(\alpha)(\Phi)]\right|_{\alpha=0}
\cdot {\mathbb C}^p,
\]
which implies that the mapping $\alpha\longmapsto {\cal B}(\alpha)$ 
is transversal (mini-transversal if (\ref{minitransv}) holds) to the
similarity orbit of $B$ at $\alpha=0$.
The conclusion now follows from Proposition \ref{ArWiprop}.
\hfill\qed

From the previous result, we are now led to define a notion 
of versal unfolding for RFDEs
of the type~(\ref{fdep}). 
\begin{Def}
The parametrized family of RFDEs~(\ref{fdep}) is said to be a 
{\em $\Lambda$-versal unfolding} (respectively a {\em $\Lambda$-mini-versal
unfolding}) for
the RFDE~(\ref{fde}) if the matrix ${\cal B}(\alpha)$ defined by~(\ref{Bunfdef}) is 
a versal unfolding (respectively a mini-versal unfolding) for $B$.
\end{Def}

\Section{Decomposition of Mat{\boldmath $_{c\times c}$} by 
{\boldmath $\Psi(0)$}\label{sec6}}

Consider now the following problem: given a linear homogeneous
RFDE
such as (\ref{fde}) and a set $\Lambda$ of solutions to 
the characteristic equation (\ref{char_eq}), 
find a $\Lambda$-versal unfolding (\ref{fdep}) for (\ref{fde}).
It is certainly not immediately obvious that this problem need admit a solution, since
the structure of the right-hand side of (\ref{cmeqsparam}) is severely
restricted by the structure of the matrix $\Psi(0)$.
To solve this problem, we will need to 
characterize the subspace of matrices
in $\mbox{\rm Mat}_{c\times c}$ whose columns are in the 
range of $\Psi(0)$, and 
show that one can build a versal unfolding of the matrix $B$ 
in (\ref{cmeqsparam})
even in
this restricted context.
We start with the following
\begin{Def}
We define ${\cal R}(\Psi(0))$ to be the set of all $c\times c$
matrices
whose columns are in the range of the matrix $\Psi(0)$.
\label{rangepsidef}
\end{Def}
The main result we will need in order to construct the above-mentioned versal unfolding 
is the following:
\begin{prop}
Let ${\cal T}:\mbox{\rm Mat}_{c\times c}\longrightarrow
\mbox{\rm Mat}_{c\times c}$ be defined by
${\cal T}(M)=[B,M]\equiv BM-MB$.
There exists a subspace $\widehat{W}\subset{\cal R}(\Psi(0))$ such that
\begin{equation}
\mbox{\rm Mat}_{c\times c}=\mbox{\rm range}({\cal T})\oplus
\widehat{\cal W}.
\label{homoleq}
\end{equation}
\label{homolprop}
\end{prop}
This entire section is devoted to proving this result, including 
giving an explicit
construction of $\widehat{\cal W}$. We will start by first
proving several lemmas, and then proceed to the proof of
Proposition \ref{homolprop}.

\subsection{Commutator of {\boldmath $B$}}
We assume that the matrix $B$ is in the following Jordan block diagonal form
\begin{equation}
B=\mbox{\rm diag}(B^1_{1}(\lambda_1),\ldots,B^1_{k_1}(\lambda_1),
B^2_1(\lambda_2),\ldots,B^2_{k_2}(\lambda_2),\ldots,B^{r}_1(\lambda_r),
\ldots,B^{r}_{k_r}(\lambda_r)),
\label{jordanB}
\end{equation}
where $\lambda_1,\ldots,\lambda_r$ are the distinct (without multiplicities) 
eigenvalues of $B$, and the Jordan block $B^j_{\ell}(\lambda_j)$ is 
$n_{j,\ell}\times n_{j,\ell}$ of the form $\lambda_j\,I_{n_{j,\ell}}+
N_{n_{j,\ell}\times n_{j,\ell}}$, where $N_{n_{j,\ell}\times n_{j,\ell}}$ 
is the nilpotent matrix with 1's on the upper diagonal, and 0's everywhere 
else. Moreover, we assume that for each $j\in \{1,\ldots,r\}$, we have
\[
n_{j,1}\geq n_{j,2}\geq\cdots\geq n_{j,k_j}.
\]
The first result we need is the following lemma, which can be found 
in \cite{Gant}:
\begin{lemma}
Let $B^*$ denote the conjugate transpose of $B$, and let $M\in
\mbox{\rm Mat}_{c\times c}$ be such that $[B^*,M]=0$.  Then $M$ 
is of the form 
\begin{equation}
M=\mbox{\rm diag}({\cal M}_1,\ldots,{\cal M}_r),
\label{Ndiagstruct}
\end{equation}
where ${\cal M}_j$ is 
$(n_{j,1}+\cdots+n_{j,k_j})\times (n_{j,1}+\cdots+n_{j,k_j})$ matrix, 
partitioned into blocks with dimensions $n_{j,p}\times n_{j,q}$, and
of the form illustrated in Figure \ref{fig1}, where each oblique segment 
in each separate block denotes a sequence of equal entries, and all 
other entries are zero.
\begin{figure}[hptb]
\centerline{%
\psfig{file=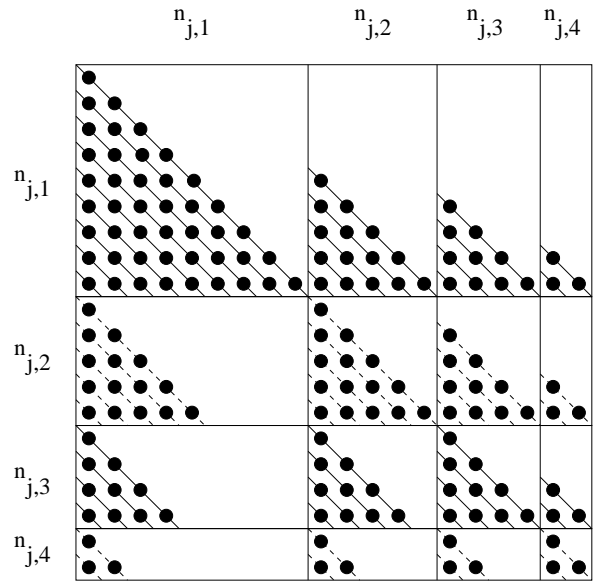,height=3.0in}}
\caption{Structure of a matrix ${\cal M}_j$.  In this example,
$k_j=4$, $n_{j,1}=9$, $n_{j,2}=5$, $n_{j,3}=4$, $n_{j,4}=2$.}
\label{fig1}
\end{figure}
\label{commstruct}
\end{lemma}
The second lemma we need is the following. The proof is left
to the appendix.
\begin{lemma}
Let $B$ be as in (\ref{jordanB}), and consider the mapping
${\cal T}:\mbox{\rm Mat}_{c\times c}\longrightarrow\mbox{\rm
  Mat}_{c\times c}$ given by ${\cal T}(M)=[B,M]$.  Then 
$Y\in\mbox{\rm range}({\cal T})$ if and only if 
$Y$ is of the form
\begin{equation}
Y=\left(\begin{array}{ccccc}
{\cal Y}_{1,1}&{\cal Y}_{1,2}&\cdots&\cdots&{\cal Y}_{1,r}\\
{\cal Y}_{2,1}&{\cal Y}_{2,2}&\cdots&\cdots&{\cal Y}_{2,r}\\
\vdots&\vdots&\vdots&\vdots&\vdots\\
\vdots&\vdots&\vdots&\vdots&\vdots\\
{\cal Y}_{r,1}&{\cal Y}_{r,2}&\cdots&\cdots&{\cal Y}_{r,r}
\end{array}
\right),
\label{rangegen}
\end{equation}
where ${\cal Y}_{p,q}$ is a $(n_{p,1}+\cdots+n_{p,k_p})\times (n_{q,1}+
\cdots+n_{q,k_q})$ matrix, with the only constraint 
being that for each $j\in\{1,\ldots,r\}$, the sub-matrix ${\cal Y}_{j,j}$
is partitioned exactly as ${\cal M}_j$ (see Figure \ref{fig1}), 
and is such that the sum of the elements in each given oblique segment 
is zero (however, in contrast to ${\cal M}_{j}$, the elements of 
${\cal Y}_{j,j}$ which are not on the oblique segments are completely 
arbitrary, i.e. not necessarily  zero).
\label{rangeTlem}
\end{lemma}
We now define useful integers which give the number of columns of 
the matrix $B$ corresponding to the first $j$ eigenvalues. 
Let $N_0=0$, and 
\[
N_{j}=N_{j-1}+\sum_{\ell=1}^{k_j}\,n_{j,\ell},\,\,\,\,
j=1,\ldots,r-1.
\]
Using Lemma~\ref{rangeTlem}, we now define a set of matrices ${\cal S}$ 
which forms a basis for $\mbox{\rm range}({\cal T})$. The set ${\cal S}$ 
consists of all $c\times c$ matrices $Y$ which are partitioned as in
(\ref{rangegen}) and whose elements are all zero except one element whose
value is 1, and possibly one other element whose value is $-1$,
satisfying the following constraints: 
\begin{itemize}
\item if the element whose value is 1 is in the block ${\cal Y}_{k,\ell}$ 
where $k\neq\ell$, then it is the {\em only} non-zero element in
the matrix $Y$,
\item if the element whose value is 1 lies in the block ${\cal Y}_{j,j}$ for
some $j\in\{1,\ldots,r\}$, then:
\begin{itemize}
\item if this element (whose value is 1) 
is not on any of the oblique segments of Figure \ref{fig1}, 
then it is the {\em only} non-zero element in
the matrix $Y$
\item if this element (whose value is 1) does lie on one of the
oblique segments of Figure \ref{fig1}, then it does not lie in any of the 
following rows
\begin{equation}\label{rows}
N_{j-1}+n_{j,1},\,\,N_{j-1}+n_{j,1}+n_{j,2},\,\,\ldots,\,\,N_{j-1}+n_{j,1}
+\cdots+n_{j,k_j},
\end{equation}
and there is a `$-1$' at the bottom end of the oblique segment which 
contains the `1'.
\end{itemize}
\end{itemize}
Note that the row numbers given by~(\ref{rows}) correspond to the last row
of $n_{j,k_j}\times n_{j,k_{j}}$ blocks as shown in Figure~\ref{fig1} 
for each $j$.

\begin{figure}[h]
\centerline{%
\psfig{file=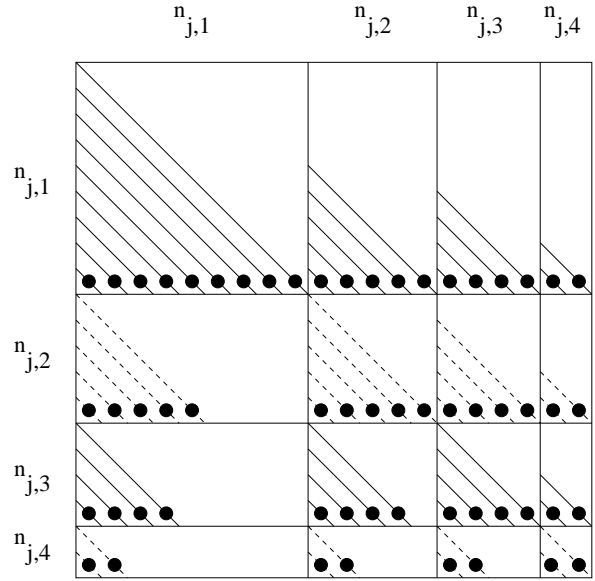,height=3.0in}}
\caption{Structure of a matrix $\omega_j$.  In this example,
$k_j=4$.  The only elements which are not forced to be zero
are those illustrated with a bold dot, at the
bottom of each oblique segment.}
\label{fig2}
\end{figure}

Now consider ${\cal W}\subset\mbox{\rm Mat}_{c\times c}$,
where ${\cal W}$ is the subspace of matrices of the form
\[\mbox{\rm diag}(\omega_1,\ldots,\omega_r),\]
where $\omega_j$ is a $(n_{j,1}+\cdots+n_{j,k_j})\times 
(n_{j,1}+\cdots+n_{j,k_j})$ matrix, partitioned into blocks with 
dimensions $n_{j,p}\times n_{j,q}$, and of the
form illustrated in Figure \ref{fig2}. In this figure, 
the only elements which are not forced to be
zero are those at the bottom of each oblique segment,
illustrated with a bold dot.  It is well-known \cite{Gant} that 
\begin{equation}
\mbox{\rm dim}({\cal W})=\sum_{j=1}^r\,\sum_{\ell=1}^{k_j}\,(2\ell-1)
n_{j,\ell}.
\label{dimunfold}
\end{equation}

It follows from the definition of the basis ${\cal S}$ of
$\mbox{\rm range}({\cal T})$ that
\begin{lemma}
\begin{equation}
\mbox{\rm Mat}_{c\times c}=\mbox{\rm range}({\cal T})\oplus {\cal W}.
\label{splitmatcc}
\end{equation}
\label{directsumTW}
\end{lemma}

\subsection{Characterizing the range of {\boldmath $\Psi(0)$}}
Before we can prove Proposition \ref{homolprop}, we need
to establish some properties about the subspace ${\cal R}(\Psi(0))$ 
defined in Definition \ref{rangepsidef}.

\begin{lemma}
The matrix $\Psi(0)=\mbox{\rm col}(\psi_1(0),\ldots,\psi_c(0))$ satisfies the
following property: for each $j\in \{1,\ldots,r\}$, 
if
$i\in \{\,N_{j-1}+n_{j,1},N_{j-1}+n_{j,1}+n_{j,2},\ldots,
N_{j-1}+n_{j,1}+\cdots+n_{j,k_j}\,\}$ then
\[
\psi_i(0)\neq 0\;\;\;\mbox{\rm and}\;\;\;\psi_i(0)\Delta(\lambda_j)=0,
\]
where the characteristic matrix $\Delta(\lambda)$ 
is as in (\ref{char_eq}).
Moreover, for each fixed $j\in\{1,\ldots,r\}$, the set of row-vectors
\[
\{\,\psi_{N_{j-1}+n_{j,1}}(0),\psi_{N_{j-1}+n_{j,1}+n_{j,2}}(0),\ldots,
\psi_{N_{j-1}+n_{j,1}+\cdots+n_{j,k_j}}(0)\,\}
\]
is linearly independent in ${\mathbb C}^{n*}$.
\label{lempsistruct}
\end{lemma}
\proof
Each row $\psi_i(\theta)$ of the matrix $\Psi(\theta)=e^{-B\theta}\Psi(0)$ 
must satisfy
\begin{equation}
\dot{\psi}_i(0)=-\int_{-\tau}^0\,\psi_i(-\theta)\,d\eta(\theta).
\label{psiconst}
\end{equation}
Since $B$ has the form (\ref{jordanB}), it follows that
for each $j\in \{1,\ldots,r\}$, if
$i\in \{\,N_{j-1}+n_{j,1},N_{j-1}+n_{j,1}+n_{j,2},\ldots,
N_{j-1}+n_{j,1}+\cdots+n_{j,k_j}\,\}$ then
the $i^{th}$ row of the matrix $e^{-B\theta}$ is zero 
except for the diagonal element which is equal to $e^{-\lambda_j\theta}$.  
It follows that for all $i\in \{\,N_{j-1}+n_{j,1},N_{j-1}+n_{j,1}+n_{j,2},
\ldots,N_{j-1}+n_{j,1}+\cdots+n_{j,k_j}\,\}$ we have
\[
\psi_i(\theta)=e^{-\lambda_j\theta}\psi_i(0),
\]
which when substituted into (\ref{psiconst}) yields 
$\psi(0)\Delta(\lambda_j)=0$. Recall that the columns of 
$\Psi(\theta)$ form a basis of $P^{*}$, therefore the set
\[
\{\,\psi_{N_{j-1}+n_{j,1}}(\theta),\psi_{N_{j-1}+n_{j,1}+n_{j,2}}(\theta),
\ldots,\psi_{N_{j-1}+n_{j,1}+\cdots+n_{j,k_j}}(\theta)\,\}
\]
is linearly independent in $C([0,\tau],{\mathbb C}^{n*})$, it follows that
\[
\{\,\psi_{N_{j-1}+n_{j,1}}(0),\psi_{N_{j-1}+n_{j,1}+n_{j,2}}(0),\ldots,
\psi_{N_{j-1}+n_{j,1}+\cdots+n_{j,k_j}}(0)\,\}
\]
is linearly independent in ${\mathbb C}^{n*}$.
\hfill
\qed
\begin{rmk}
For each $j\in\{1,\ldots,r\}$, let 
\[
\{\,\psi_{N_{j-1}+n_{j,1}}(0),\psi_{N_{j-1}+n_{j,1}+n_{j,2}}(0),\ldots,
\psi_{N_{j-1}+n_{j,1}+\cdots+n_{j,k_j}}(0)\,\}
\]
be as in Lemma \ref{lempsistruct}.  Consider the linear mapping
\[
\Pi_j:{\mathbb C}^{n}\longrightarrow {\mathbb C}^{k_j}
\]
defined by
\[
\Pi_j(v)=\mbox{\rm col}
(\psi_{N_{j-1}+n_{j,1}}(0),\psi_{N_{j-1}+n_{j,1}+n_{j,2}}(0),\ldots,
\psi_{N_{j-1}+n_{j,1}+\cdots+n_{j,k_j}}(0))\cdot v.
\]
Then it follows from Lemma \ref{lempsistruct} that $\Pi_j$ is onto ${\mathbb C}^{k_j}$.
\label{rankpsilem}
\end{rmk}

\subsection{Proof of Proposition \ref{homolprop}}
We are now ready to prove Proposition \ref{homolprop}.  
For purposes of clarity, we first prove it in the case where $r=1$ in
(\ref{jordanB}), and then we show how to generalize the arguments of this
proof to the case $r>1$. Before giving the proof, we first 
establish some useful notation.

In the case where $r=1$ in (\ref{jordanB}), any matrix 
$M\in\mbox{\rm Mat}_{c\times c}$ can be partitioned as
\begin{equation}
M=\left(\begin{array}{ccccc}
{\cal M}_{1,1}&{\cal M}_{1,2}&\cdots&\cdots&{\cal
  M}_{1,k_1}\\
{\cal M}_{2,1}&{\cal M}_{2,2}&\cdots&\cdots&{\cal
  M}_{2,k_1}\\
\vdots&\vdots&\vdots&\vdots&\vdots\\
\vdots&\vdots&\vdots&\vdots&\vdots\\
{\cal M}_{k_1,1}&{\cal M}_{k_1,2}&\cdots&\cdots&{\cal
  M}_{k_1,k_1}
\end{array}
\right),
\label{Mdecompj1}
\end{equation}
where ${\cal M}_{\xi,\lambda}$ is $n_{1,\xi}\times n_{1,\lambda}$, 
for $\xi,\lambda\in\{1,\ldots,k_1\}$. It is convenient to label the 
elements of $M$ according to this partitioning as follows: for 
$\xi,\lambda\in\{1,\ldots,k_1\}$, $m\in\{1,\ldots,n_{1,\lambda}\}$ and 
$u\in\{1,\ldots,n_{1,\xi}\}$, we denote by $M^{\xi,\lambda,u,m}$ the 
element of the matrix $M$ which lies in the block ${\cal M}_{\xi,\lambda}$ 
in (\ref{Mdecompj1}), at the intersection of row $u$ and column 
$n_{1,\lambda}-m+1$ relative to this block.

Fix $\xi,\lambda\in\{1,\ldots,k_1\}$ and let ${\cal Q}(\xi,\lambda)$ be
the set of integers $m\in\{1,\ldots,n_{1,\lambda}\}$ such that
column $n_{1,\lambda}-m+1$ intersects an oblique line in the block
${\cal M}_{\xi,\lambda}$. In particular,
\[
{\cal Q}(\xi,\lambda)=\left\{\begin{array}{ll}
\{1,\ldots,n_{1,\lambda}\} & \mbox{if\,} n_{1,\xi}\geq n_{1,\lambda}\\
\{n_{1,\lambda}-n_{1,\xi}+1,\ldots,n_{1,\lambda}\} & \mbox{if\,} n_{1,\xi}< n_{1,\lambda}.
\end{array}\right.
\]
Now, for all $\xi,\lambda \in\{1,\ldots,k_1\}$ and 
$m\in {\cal Q}(\xi,\lambda)$, we define $\Omega_{\xi,\lambda,m}$
as the $c\times c$ matrix whose entries are all zero except
$\Omega_{\xi,\lambda,m}^{\xi,\lambda,n_{1,\xi},m}=1$.  It follows that
the set of all matrices $\Omega_{\xi,\lambda,m}$ thus defined is a
basis for the subspace ${\cal W}$ in (\ref{splitmatcc}) 
(see Figure \ref{fig2}).

Similarly, fix $\xi,\lambda \in\{1,\ldots,k_1\}$ and 
$m\in\{1,\ldots,n_{1,\lambda}\}$. Let ${\cal P}(\xi,\lambda,m)$
consist of all $u\in\{1,\ldots,n_{1,\xi}-1\}$ such that row $u$ has
an intersection with an oblique line in the block $(\xi,\lambda)$ at
column $n_{1,\lambda}-m+1$. 
Of course, if $n_{1,\xi}=1$, $m\not\in {\cal Q}(\xi,\lambda)$ or 
$m=\min {\cal Q}(\xi,\lambda)$ then ${\cal P}(\xi,\lambda,m)=\emptyset$.

\begin{figure}\label{fig:matnot}
\psfrag{m}{$m=6$}
\psfrag{P}{${\cal P}(\xi,\lambda,6)$}
\psfrag{Q}{${\cal Q}(\xi,\lambda)$}
\psfrag{p42}{$p(6,3)$}
\centerline{%
\epsfxsize=0.6\hsize\epsffile{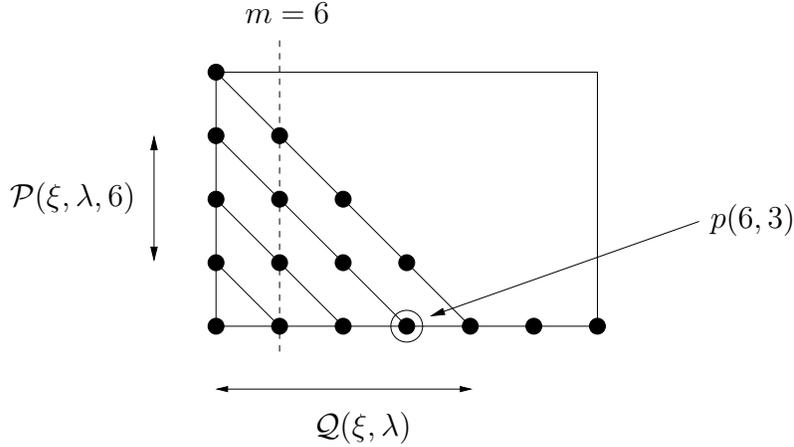}}
\caption{Example of sets ${\cal Q}$, ${\cal P}$ and of a point $p(m,u)$
for a $5\times 7$ matrix. Here ${\cal Q}=\{3,4,5,6,7\}$ and 
${\cal P}=\{2,3,4\}$.}
\end{figure}

For all $\xi,\lambda\in\{1,\ldots,k_1\}$, choose
$m\in {\cal Q}(\xi,\lambda)$ such that ${\cal P}(\xi,\lambda,m)$ is nonempty
and choose $u\in {\cal P}(\xi,\lambda,m)$. 
Define the $c\times c$ matrix $E_{\xi,\lambda,u,m}$ with only two nonzero 
entries as follows. Let $E_{\xi,\lambda,u,m}^{\xi,\lambda,u,m}=1$, and
$E_{\xi,\lambda,u,m}^{\xi,\lambda,n_{1,\xi},p(m,u)}=-1$,
where $p(m,u)=m-n_{1,\xi}+u$. Note that $p(m,u)$ is always less than $m$.
See Figure~\ref{fig:matnot} for an illustration of ${\cal Q}(\xi,\lambda)$, 
${\cal P}(\xi,\lambda,m)$ and $p(m,u)$.
Each of these matrices $E_{\xi,\lambda,u,m}$ defined above is in 
$\mbox{\rm range}({\cal T})$ (Lemma \ref{rangeTlem}).
\begin{rmk}
It follows from the definitions of the matrices 
$\Omega_{\xi,\lambda,m}$ and
the matrices $E_{\xi,\lambda,u,m}$ that
if $\lambda\in\{1,\ldots,k_1\}$, $\xi\in\{1,\ldots,k_1\}$,
$m\in\{1,\ldots,n_{1,\lambda}\}$ and
$u\in {\cal P}(\xi,\lambda,m)$, then $p(m,u)\in {\cal Q}(\xi,\lambda)$ and
$E_{\xi,\lambda,u,m}+\Omega_{\xi,\lambda,p(m,u)}$ is a
matrix whose only non-zero element is
\[
(E_{\xi,\lambda,u,m}+\Omega_{\xi,\lambda,p(m,u)})^{\xi,\lambda,u,m}=1.
\]
\label{bascomb}
\end{rmk}

We now proceed to the proof of Proposition \ref{homolprop} in the
case where $r=1$.
\begin{lemma}
Suppose $B$ in (\ref{jordanB}) is such that $r=1$.
Then Proposition \ref{homolprop} holds.
\label{homolproplem2}
\end{lemma}
\proof
We begin by defining a subspace of matrices contained in ${\cal R}(\Psi(0))$
and isomorphic to ${\cal W}$. We do this in the following way.

By virtue of Remark \ref{rankpsilem}, let $v_1,\ldots,v_{k_1}\in 
{\mathbb C}^n$ be such that $\Pi_1(v_{\ell})$ is a $k_1$-dimensional 
vector whose only non-zero component is the $\ell^{th}$ 
component, whose value is 1. Now, consider the $c$-dimensional vector 
$\Psi(0)v_{\ell}$.  We label the components of this vector using the 
integers $\xi\in\{1,\ldots,k_1\}$ and $u\in\{1,\ldots,n_{1,\xi}\}$ as 
follows: $(\Psi(0)v_{\ell})^{\xi,u}$ is the $u^{th}$ component
of $\Psi(0)v_{\ell}$ if $\xi=1$; and the 
$(n_{1,1}+\cdots +n_{1,\xi-1}+u)^{th}$ component of 
$\Psi(0)v_{\ell}$ if $\xi>1$. We then have
\begin{equation}
(\Psi(0)\,v_{\ell})^{\xi,u}=\left\{\begin{array}{cl}
\gamma_{\xi,u,\ell}&\mbox{\rm if}\,\,u\notin \{n_{1,1},\ldots,n_{1,k_1}\}\\[0.2in]
1&\mbox{\rm if}\,\,\xi=\ell\,\,\mbox{\rm and}\,\,u=n_{1,\xi}\\[0.2in]
0&\mbox{\rm otherwise}
\end{array}\right.
\label{psiv2}
\end{equation}
where the exact values of the coefficients $\gamma_{\xi,u,\ell}$ are not 
important.

Now for all $\xi,\lambda\in\{1,\ldots,k_1\}$ and $m\in{\cal Q}(\xi,\lambda)$ 
define $R_{\xi,\lambda,m}$ to be the $n\times c$ matrix whose 
$(n_{1,1}+\cdots+n_{1,\lambda}-m+1)^{th}$ column is $v_{\xi}$, and all 
other columns are zero.
We then define a linear mapping ${\cal E}:{\cal W}\longrightarrow
{\cal R}(\Psi(0))$ by the following action on the basis elements of
${\cal W}$:
\begin{equation}
{\cal E}(\Omega_{\xi,\lambda,m})=\Psi(0)
R_{\xi,\lambda,m},\;\;\;\;
\xi,\lambda\in\{1,\ldots,k_1\},\,m\in {\cal Q}(\xi,\lambda).
\label{defE}
\end{equation}
It is clear that ${\cal E}$ is an isomorphism between ${\cal W}$ and
$\widehat{\cal W}\equiv {\cal E}({\cal W})\subset {\cal R}(\Psi(0))$,
since the set
\[
\{\,{\cal E}(\Omega_{\xi,\lambda,m})\,\,:\,\,
\xi,\lambda\in\{1,\ldots,k_1\},\,m\in {\cal Q}(\xi,\lambda)\,\}
\]
is linearly independent in $\mbox{\rm Mat}_{c\times c}$.

Our strategy now is to show that
\[
{\cal W}\subset\mbox{\rm range}({\cal T})+\widehat{\cal W},
\]
from which it follows from (\ref{splitmatcc}) that
\[
\mbox{\rm Mat}_{c\times c}=\mbox{\rm range}({\cal T})+\widehat{\cal W},
\]
and since ${\cal W}$ and $\widehat{\cal W}$ are isomorphic, we will then have
\begin{equation}
\mbox{\rm Mat}_{c\times c}=\mbox{\rm range}({\cal T})\oplus\widehat{\cal W}.
\label{splitmatcc3}
\end{equation}

\begin{figure}[hptb]
\centerline{%
\psfig{file=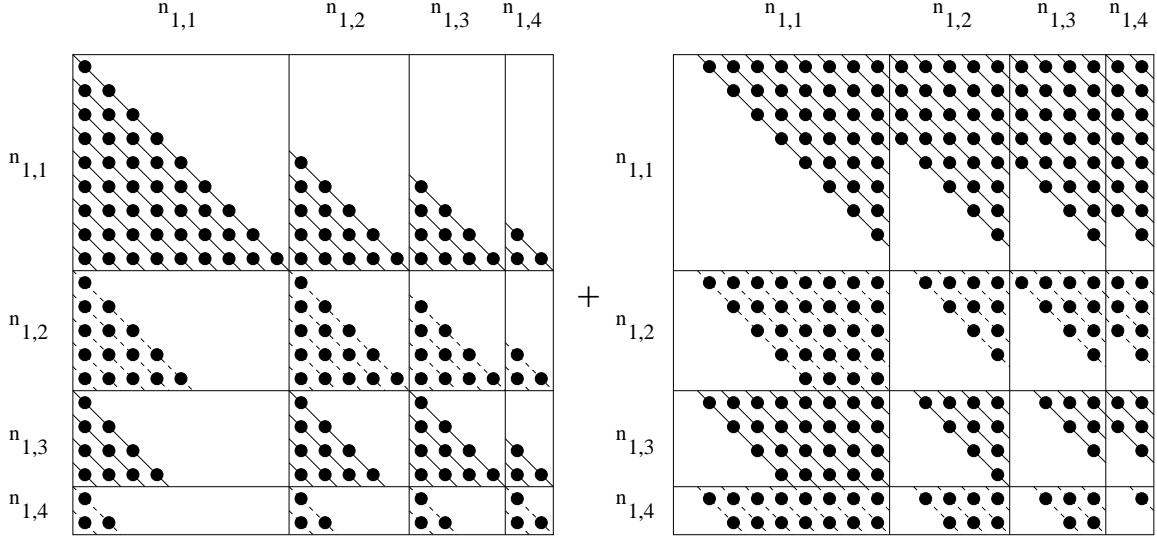,width=6.0in}}
\caption{Decomposition of $Z\in\mbox{\rm Mat}_{c\times c}$.  In this example,
$k_1=4$.  The elements which are not on the oblique segments are all
zero. The first summand is $\Gamma(Z)$.
The second summand is in $\mbox{\rm range}({\cal T})$.}
\label{fig3}
\end{figure}

Any matrix $Z\in\mbox{\rm Mat}_{c\times c}$ can be written as
a sum of two matrices as is illustrated in Figure \ref{fig3}, where the
elements which are not on the oblique segments are all zero. It is clear 
from Lemma \ref{rangeTlem} that the second summand in Figure \ref{fig3} 
belongs to $\mbox{\rm range}({\cal T})$.  Thus, for any 
$Z\in\mbox{\rm Mat}_{c\times c}$, we define $\Gamma(Z)$ to be the first 
summand in this decomposition, as illustrated in Figure \ref{fig3}.  
It immediately follows that for any $Z\in\mbox{\rm Mat}_{c\times c}$,
we have $Z-\Gamma(Z)\in\mbox{\rm range}({\cal T})$.

Fix a value of $\lambda\in\{1,\ldots,k_1\}$.  
For all $\xi\in\{1,\ldots,k_1\}$, if $1\in {\cal Q}(\xi,\lambda)$
then $\Gamma(\Psi(0) R_{\xi,\lambda,1})=
\Omega_{\xi,\lambda,1}$, from which it follows that
\[
\Omega_{\xi,\lambda,1}\in\mbox{\rm range}({\cal T})+\widehat{\cal
  W},\;\;\;
\forall\,\xi\in\{1,\ldots,k_1\}\,\,\mbox{\rm such that}\,1\in {\cal
  Q}(\xi,\lambda).
\]
If $n_{1,\lambda}>1$, let $\mu$ be an integer in the
range $1,\ldots,n_{1,\lambda}-1$ such that for all $m=1,\ldots,\mu$,
we have
\begin{equation}
\Omega_{\xi,\lambda,m}\in\mbox{\rm range}({\cal T})+\widehat{\cal
  W},\;\;\;
\forall\,\xi\in\{1,\ldots,k_1\}\,\,\mbox{\rm such that}\,m\in {\cal
  Q}(\xi,\lambda).
\label{indhyp2}
\end{equation}
We claim that this implies that 
\begin{equation}
\Omega_{\xi,\lambda,\mu+1}\in\mbox{\rm range}({\cal T})+\widehat{\cal
  W},\;\;\;\forall\,\xi\in\{1,\ldots,k_1\}\,\,\mbox{\rm such that}\,
\mu+1\in {\cal Q}(\xi,\lambda).
\label{indconc2}
\end{equation}
From this claim, and the fact that the above arguments are independent
of the particular choice of $\lambda$, it follows by induction that 
(\ref{splitmatcc3}) holds. It thus remains to prove the claim.

Suppose that $\xi\in\{1,\ldots,k_1\}$ is such that
$\mu+1\in {\cal Q}(\xi,\lambda)$.  Then a simple computation shows that
\begin{equation}
\Gamma(\Psi(0)R_{\xi,\lambda,\mu+1})=\Omega_{\xi,\lambda,\mu+1}+
G_{\xi,\lambda,\mu+1},
\label{psiog}
\end{equation}
where $G_{\xi,\lambda,\mu+1}$ is a $c\times c$ matrix whose columns
are all zero except possibly the $(n_{1,1}+\cdots +n_{1,\lambda}
-\mu)^{th}$ column, whose $(\chi,u)$ component is given by the formula
\[
\left\{
\begin{array}{cl}
\gamma_{\chi,u,\xi}&\mbox{\rm if}\,\,u\in {\cal P}(\chi,\lambda,\mu+1)\\
[0.2in] 0&\mbox{\rm otherwise},
\end{array}
\right.
\]
where the coefficients $\gamma_{\chi,u,\xi}$ are as in
(\ref{psiv2}). Now, $\Psi(0)R_{\xi,\lambda,\mu+1}
-\Gamma(\Psi(0)R_{\xi,\lambda,\mu+1})\in \mbox{\rm range}({\cal T})$
implies that $\Gamma(\Psi(0)R_{\xi,\lambda,\mu+1})\in 
\mbox{\rm range}({\cal T}) + \widehat{\cal W}$. 
From Remark~\ref{bascomb},
\begin{equation}
G_{\xi,\lambda,\mu+1}=\sum_{\chi=1}^{k_1}\,\sum_{u\,\in\, 
{\cal P}(\chi,\lambda,\mu+1)}\,\gamma_{\chi,u,\xi}\,
(\,E_{\chi,\lambda,u,\mu+1}+\Omega_{\chi,\lambda,p(\mu+1,u)}\,).
\label{Gdec}
\end{equation}
from which it follows that (\ref{indconc2}) holds since
$E_{\chi,\lambda,u,\mu+1}\in \mbox{\rm range}({\cal T})$
and $\Omega_{\chi,\lambda,p(\mu+1,u)} \in \mbox{\rm range}({\cal T})
+\widehat{\cal W}$ from the fact that $p(\mu+1,u)<\mu+1$ and by the 
induction hypothesis. 
\hfill
\qed

\noindent
{\bf Proof of Proposition \ref{homolprop}}\hspace*{0.2in}
The essential idea here is to decompose the proof into $r$ separate
blocks, where we use the arguments of the proof of Lemma \ref{homolproplem2} 
on each of the blocks.

First, for all $j\in\{1,\ldots,r\}$, we let ${\cal W}_j$
denote the subspace of ${\cal W}$ consisting of matrices
whose columns are all zero except the columns between
$N_{j-1}+1$ and $N_j$ inclusively.

As in the proof of Lemma \ref{homolproplem2}, we construct
a basis $\{\,\Omega_{j;\xi,\lambda,u,m}\,\}$ of ${\cal W}_j$,
and we note that the space ${\cal W}$ in (\ref{splitmatcc}) is equal to
${\cal W}_1\oplus\cdots\oplus {\cal W}_r$.
From Remark \ref{rankpsilem}, for each $j\in\{1,\ldots,r\}$ we can
choose vectors $v_{j,1},\ldots,v_{j,k_j}$ which are such that
the vector $\Pi_j(v_{j,\ell})$ has a 1 in row $\ell$ and 0's
everywhere else.  
For each $j\in\{1,\ldots,r\}$, we then construct matrices 
$R_{j;\xi,\lambda,m}$ in a similar manner as we did in
Lemma \ref{homolproplem2}, whose only non-zero column
is equal to one of the vectors $v_{j,1},\ldots,v_{j,k_j}$
described above.  We then define linear mappings ${\cal E}_j:{\cal W}_j
\longrightarrow {\cal R}(\Psi(0))$ as in (\ref{defE}), and it follows
that ${\cal W}_j$ is isomorphic to 
$\widehat{\cal W}_j={\cal E}_j({\cal W}_j)$.

Now, any $Z\in\mbox{\rm Mat}_{c\times c}$ can be partitioned as in
(\ref{rangegen}), i.e.
\begin{equation}
Z=\left(\begin{array}{ccccc}
{\cal Z}_{1,1}&{\cal Z}_{1,2}&\cdots&\cdots&{\cal Z}_{1,r}\\
{\cal Z}_{2,1}&{\cal Z}_{2,2}&\cdots&\cdots&{\cal Z}_{2,r}\\
\vdots&\vdots&\vdots&\vdots&\vdots\\
\vdots&\vdots&\vdots&\vdots&\vdots\\
{\cal Z}_{r,1}&{\cal Z}_{r,2}&\cdots&\cdots&{\cal Z}_{r,r}
\end{array}
\right),
\label{genZ}
\end{equation}
where ${\cal Z}_{p,q}$ 
is a $(n_{p,1}+\cdots+n_{p,k_p})\times (n_{q,1}+\cdots+n_{q,k_q})$ matrix.
Rewrite $Z$ as
\begin{equation}
Z=\left(\begin{array}{ccccc}
{\cal Z}_{1,1}&0&\cdots&\cdots&0\\
0&{\cal Z}_{2,2}&\cdots&\cdots&0\\
\vdots&\vdots&\vdots&\vdots&\vdots\\
\vdots&\vdots&\vdots&\vdots&\vdots\\
0&0&\cdots&\cdots&{\cal Z}_{r,r}
\end{array}
\right)+
\left(\begin{array}{ccccc}
0&{\cal Z}_{1,2}&\cdots&\cdots&{\cal Z}_{1,r}\\
{\cal Z}_{2,1}&0&\cdots&\cdots&{\cal Z}_{2,r}\\
\vdots&\vdots&\vdots&\vdots&\vdots\\
\vdots&\vdots&\vdots&\vdots&\vdots\\
{\cal Z}_{r,1}&{\cal Z}_{r,2}&\cdots&\cdots&0
\end{array}
\right),
\label{Zdecomp}
\end{equation}
where the second summand on the right-hand side of (\ref{Zdecomp}) belongs
to $\mbox{\rm range}({\cal T})$ (Lemma \ref{rangeTlem}).
Furthermore, following the same method as in Lemma \ref{homolproplem2},
for each $j\in\{1,\ldots,r\}$ we write ${\cal Z}_{j,j}$
as a sum of matrices as is illustrated in Figure \ref{fig3},
where the second summand belongs to 
$\mbox{\rm range}({\cal T})$.  Thus, for $Z$ as in (\ref{genZ}), we
define $\Gamma(Z)$ as the first summand on the right-hand side of
(\ref{Zdecomp}), where the blocks ${\cal Z}_{j,j}$ are of the form of the
first summand in Figure \ref{fig3}.  
It follows that for any $Z\in\mbox{\rm Mat}_{c\times c}$,
we have $Z-\Gamma(Z)\in\mbox{\rm range}({\cal T})$.

Finally, we use an induction argument similar to that used in the proof
of Lemma \ref{homolproplem2} to show that for each $j\in\{1,\ldots,r\}$,
we have ${\cal W}_j\subset\mbox{\rm range}({\cal T})+\widehat{\cal W}_j$.  
We thus define $\widehat{\cal W}=\widehat{\cal W}_1\oplus
\cdots\oplus\widehat{\cal W}_r$. This completes the proof.
\hfill
\qed

\Section{Construction of a {\boldmath $\Lambda$}-versal unfolding}

In this section, we will solve the problem which was posed at the
beginning of Section \ref{sec6}; that is,
given a linear homogeneous
RFDE
such as (\ref{fde}) and a set $\Lambda$ of solutions to 
the characteristic equation (\ref{char_eq}), 
find a $\Lambda$-versal unfolding (\ref{fdep}) for (\ref{fde}).

\subsection{Main result}

The following useful result is proven in \cite{Hal79} and \cite{FMR}.
\begin{lemma}
Let $I$ be an interval in ${\mathbb R}$ and suppose that $\varphi_1,
\ldots,\varphi_c \in C(I,{\mathbb C}^n)$ are linearly independent 
functions.  For each $\theta\in I$, define
$\Phi_c(\theta)=(\varphi_1(\theta),\ldots,\varphi_c(\theta))$.  
For a fixed $\tau_0\in I$, denote by $q$ the rank of the $n\times c$ 
matrix $\Phi_c(\tau_0)$. Then there exist $c-q$
distinct points $\tau_1,\ldots,\tau_{c-q}\in I\setminus\{\tau_0\}$ 
such that the $n(c-q+1)\times c$ matrix 
$\mbox{\rm col}(\Phi_c(\tau_0),\ldots,\Phi_c(\tau_{c-q}))$ has rank $c$.
\label{lemFM}
\end{lemma}

Now, suppose $\varphi_1,\ldots,\varphi_c$ are the basis elements of the space
$P$ in (\ref{splitcn}), and let
$\Phi(\theta)=(\varphi_1(\theta),\ldots,\varphi_c(\theta))$, and $I=[-\tau,0]$.
Since $\dot{\Phi}=\Phi B$, it follows that 
the rows of $\Phi$ are solutions to the
vector ordinary differential equation $\dot{\xi}=\xi\,B$, and so the rank of 
the matrix $\Phi(\theta)$ is independent of $\theta$ \cite{Har}.
We now have the following
\begin{prop}
Let $q=\mbox{\rm rank}(\Phi(0))$, and let $\tau_0,\ldots,\tau_{c-q}$ be as in 
Lemma \ref{lemFM}.  Then for any $n\times c$ matrix $R$, there exist
matrices $A_0,\ldots,A_{c-q}\in\mbox{\rm Mat}_{n\times n}$ such that
\begin{equation}
R=\sum_{j=0}^{c-q}\,A_j\Phi(\tau_j).
\label{unfcand}
\end{equation}
\label{unfcandprop}
\end{prop}
\proof
Define the linear mapping ${\cal K}:\mbox{\rm Mat}_{n\times (n(c-q+1))}
\longrightarrow
\mbox{\rm Mat}_{n\times c}$ by
\[
{\cal K}({\cal A})={\cal A}\cdot\mbox{\rm col}(\Phi(\tau_0),\ldots,
\Phi(\tau_{c-q})).
\]
If we associate 
$\mbox{\rm Mat}_{n\times (n(c-q+1))}\cong {\mathbb C}^{n(n(c-q+1))}$
and $\mbox{\rm Mat}_{n\times c}\cong {\mathbb C}^{nc}$ in the standard way, 
then the $nc\times n(n(c-q+1))$ matrix representation of ${\cal K}$ is given by
\[
{\cal K}\sim
I_{n}\otimes (\mbox{\rm col}(\Phi(\tau_0),\ldots,\Phi(\tau_{c-q})))^T, 
\]
whose rank is $nc$.
Thus, ${\cal K}$ is onto $\mbox{\rm Mat}_{n\times c}$.  For a given 
$R\in \mbox{\rm Mat}_{n\times c}$, let ${\cal A}\in 
\mbox{\rm Mat}_{n\times (n(c-q+1))}$ be such that ${\cal K}({\cal A})=R$.  The
conclusion of the Proposition follows by partitioning the
matrix ${\cal A}$ as ${\cal A}=(A_0,\ldots,A_{c-q})$, where the
$A_j$ are $n\times n$ matrices, $j=0,\ldots,c-q$.
\hfill
\qed

\begin{rmk}
The crucial element 
in the proof of Proposition \ref{unfcandprop} is the fact that the
rank of the matrix
$\mbox{\rm col}(\Phi(\tau_0),\ldots,\Phi(\tau_{c-q}))$ is equal to
$c$.  Lemma \ref{lemFM} assures us that we can always achieve this
if we use $q=\mbox{\rm rank}(\Phi(0))$ and if the
delay times $\tau_0,\ldots,\tau_{c-q}$ are chosen appropriately.  
However, as is noted in
\cite{FMR}, in certain cases it is possible to achieve 
$\mbox{\rm rank}(\mbox{\rm
  col}(\Phi(\tau_0),\ldots,\Phi(\tau_{c-q})))=c$
with a value of $q$ larger than $\mbox{\rm rank}(\Phi(0))$.
For example, it is possible to have the linear RFDE (\ref{fde})
on ${\mathbb C}^2$ such that $P$ is 4-dimensional with basis matrix
\[
\Phi(\theta)=\left(\begin{array}{cccc}
e^{i\omega_1\theta}&e^{-i\omega_1\theta}&0&0\\
0&0&e^{i\omega_2\theta}&e^{-i\omega_2\theta}
\end{array}
\right),
\]
where $\omega_1$ and $\omega_2$ are distinct real numbers.
Thus, we have $c=4$ and $\mbox{\rm rank}(\Phi(0))=2$.  However, for
any
$\tau_0$ and $\tau_1$ such that
$\tau_0\neq\tau_1$, we have
\[
\mbox{\rm rank}(\mbox{\rm col}(\Phi(\tau_0),\Phi(\tau_1)))=
\mbox{\rm rank}\left(\begin{array}{cccc}
e^{i\omega_1\tau_0}&e^{-i\omega_1\tau_0}&0&0\\
0&0&e^{i\omega_2\tau_0}&e^{-i\omega_2\tau_0}\\
e^{i\omega_1\tau_1}&e^{-i\omega_1\tau_1}&0&0\\
0&0&e^{i\omega_2\tau_1}&e^{-i\omega_2\tau_1}
\end{array}
\right)=4.
\]
\label{numdels}
\end{rmk}

From Proposition \ref{homolprop}, we know that
\[
\mbox{\rm Mat}_{c\times c}=\mbox{\rm range}({\cal T}) \oplus
\widehat{\cal W},
\]
where $\widehat{\cal W}$ is isomorphic to ${\cal W}$ (see
(\ref{splitmatcc})), and the isomorphism ${\cal E}$ is described
in the proof of Lemma \ref{homolproplem2} and
Proposition \ref{homolprop}.
Let $\{\Omega_1,\ldots\Omega_{\delta}\}$ be a basis for ${\cal W}$.
From Proposition \ref{unfcandprop}, there exist matrices
$A^m_0,\ldots,A^m_{c-q}\in\mbox{\rm Mat}_{n\times n}$ such that
\begin{equation}
{\cal E}(\Omega_m)=\sum_{j=0}^{c-q}\,\Psi(0)\,A^m_j\,\Phi(\tau_j),
\;\;\;m=1,\ldots,\delta,
\label{eqforAs}
\end{equation}
where
the set $\{\,{\cal E}(\Omega_m)\,\}$
spans a complement to
$\mbox{\rm range}({\cal T})$ in $\mbox{\rm Mat}_{c\times c}$.
Clearly the previous statement still holds if instead of being 
as previously defined,
${\cal E}$ is
any injective linear mapping from ${\cal W}$ into ${\cal
  R}(\Psi(0))$
whose range is a complement to $\mbox{\rm range}({\cal T})$ in
$\mbox{\rm Mat}_{c\times c}$.

The following is our main result of this section, and follows from
Theorem \ref{thmsuff}, Propositions \ref{homolprop}
and \ref{unfcandprop}, and the previous discussion.
\begin{thm}
Let $\{\Omega_1,\ldots,\Omega_{\delta}\}$ be a basis for ${\cal W}$
(see (\ref{splitmatcc})), and let $\tau_0,\ldots,\tau_{c-q}$ be such that
\[
\mbox{\rm rank}(\mbox{\rm
  col}(\Phi(\tau_0),\ldots,\Phi(\tau_{c-q})))=c.
\]
Let ${\cal E}$ be an injective linear mapping from ${\cal W}$ into
${\cal R}(\Psi(0))$, such that 
\[
\mbox{\rm Mat}_{c\times c}=\mbox{\rm range}({\cal T})\oplus
{\cal E}({\cal W}),
\]
(there exists at least one such ${\cal E}$),
 and let $A^m_j$ ($m\in\{1,\ldots,\delta\}$,
$j\in\{0,\ldots,c-q\}$) be $n\times n$ matrices which solve
(\ref{eqforAs}).  For each $m\in\{1,\ldots,\delta\}$, let
${L}_m$ be the bounded linear operator from
$C([-\tau,0],{\mathbb C}^n)$ into ${\mathbb C}^n$ defined by
\begin{equation}
{L}_m(z)=\sum_{j=0}^{c-q}\,A^m_j\,z(\tau_j).
\label{Lmdef}
\end{equation}
Let ${\cal L}(\alpha)$ be the $\delta$-parameter family of
bounded linear operators from
$C([-\tau,0],{\mathbb C}^n)$ into ${\mathbb C}^n$ defined by
\begin{equation}
{\cal L}(\alpha)={\cal L}_0+\sum_{m=1}^{\delta}\,\alpha_m\,L_m
\label{Lalphaconst}
\end{equation}
where the $\alpha_m$ are complex parameters, and
${\cal L}_0$ is as in (\ref{fde}). 
Then (\ref{fdep}) is a $\Lambda$-mini-versal unfolding of (\ref{fde}).
\label{versunfconstthm}
\end{thm}

\subsection{First order scalar equations}
In the case of first order scalar linear RFDEs, that is
$z_t\in C_{1}$ in (\ref{fde}), we are interested in the
following question. Rewrite~(\ref{Lalphaconst}) as follows
\begin{equation}
{\cal L}(\alpha)={\cal L}_0+\sum_{j=0}^{c-q} 
\left(\sum_{m=1}^{\delta}\,\alpha_m A_{j}^{m}\right)z(\tau_{j}).
\label{Lformula}
\end{equation}
We wish to 
show that it is possible to
find a change of coordinates which simplifies~(\ref{Lformula}) to
\begin{equation}\label{simpLform}
{\cal L}(\beta)={\cal L}_0+\sum_{j=0}^{c-q} \beta_{j} z(\tau_{j})
\end{equation}
where $\beta_{j}\in {\mathbb C}$ for all $j=0,\ldots,c-q$.
Recall that $B$ is a $c\times c$ matrix and $\delta=\dim({\cal W})$. 
We have the following result.
\begin{prop}
If (\ref{fde}) is a scalar equation, then $\delta=c$.
\end{prop}

\proof By Theorem 5.1 of~\cite{FMR}, since $n=1$, the number of Jordan
blocks for each eigenvalue of $B$ is $1$. Thus $c=n_{1,1}+\cdots+n_{r,1}$. 
Now, $\delta$ is given by~(\ref{dimunfold}) with $k_{j}=1$ for all 
$j=1,\ldots,r$. Hence the equality holds.\qed

\begin{thm}
Consider a $\Lambda$-versal unfolding of a scalar 
equation~(\ref{fde}) given by~(\ref{Lalphaconst}). Then there exists 
a change of coordinates in parameter space ${\mathbb C}^{\delta}$ which 
brings (\ref{Lformula}) to~(\ref{simpLform}) where $q=1$, 
$\beta=(\beta_0,\ldots,\beta_{c-1})$ with $\beta_{j}\in {\mathbb C}$ 
for $j=0,\ldots,c-1$. 
\label{thm:scalar}
\end{thm}
\proof Since $n=1$, $q=\mbox{rank}(\Phi(0))=1$ and therefore $c$ delays
$\tau_0,\ldots,\tau_{c-1}$ are necessary to solve equation~(\ref{unfcand}).
Let $R_1,\ldots,R_\delta$ be the $1\times c$ matrices such that 
${\cal E}(\Omega_m)=\Psi(0)R_{m}$. By Proposition~\ref{unfcandprop},
\[
R_{m}=\sum_{j=0}^{c-1} A_{j}^{m} \Phi(\tau_{j})
\]
with $A_{j}^{m}\in {\mathbb C}$. We can rewrite this equation
as a matrix equation
\[
R_{m}^{T}=(\Phi(\tau_0)^{T},\ldots,\Phi(\tau_{c-1})^{T})A^{m}\equiv 
{\cal P}A^{m}.
\]
where $A^{m}=(A_{0}^{m},\ldots,A_{c-1}^{m})^{T}$ and $^T$ is transposition. 
By choice of $\tau_0,\ldots,\tau_{c-1}$, the determinant of ${\cal P}$ is 
nonzero so that $A^{m}={\cal P}^{-1}R_{m}^{T}$.

Let $\beta=(\beta_0,\ldots,\beta_{c-1})^{T}$ and 
$\alpha=(\alpha_1,\ldots,\alpha_{c})^{T}$. Since $\delta=c$, 
set
\[
\beta=\left(A^{1},\ldots,A^{c}\right)\alpha.
\]
Now $\left(A^{1},\ldots,A^{c}\right)={\cal P}^{-1}\left(R_1^T,\ldots,
R_c^T\right)$ is nonsingular since ${\cal P}^{-1}$ is a nonsingular
$c\times c$ matrix and we claim that $\left(R_1^T,\ldots,R_c^T\right)$ is
also a nonsingular $c\times c$ matrix. Hence, this change of coordinates 
yields the result.

We now prove the claim. Again by Theorem 5.1 of~\cite{FMR}, the number of 
Jordan blocks for each eigenvalue of $B$ is $1$. The set of $1\times c$
vectors $\{R_1,\ldots,R_{\delta}\}$ is partitioned into subsets 
$\{R_{j;1,1,m}|\,m\in {\cal Q}(1,1)\}$, for each $j\in\{1,\ldots,r\}$ and 
defined as in the proof of Lemma~\ref{homolproplem2} and 
Proposition~\ref{homolprop}. That is, the only nonzero element in 
$R_{j;1,1,m}$ is in the $(n_{1,1}+n_{2,1}+\ldots+n_{j,1}-m+1)^{th}$ column. 
Hence the vectors $R_1,\ldots,R_\delta$ are linearly independent since
the unique nonzero element for each vector lies in a different column.
\qed

\subsection{Independence of unfolding on choice of basis}

Suppose we perform changes of bases in the spaces $P$ and $P^*$ (see (\ref{splitcn})),
$\Psi^{\#}=U^{-1}\Psi$ and $\Phi^{\#}=\Phi\,U$, where $U$ is an invertible
$c\times c$ matrix.  Then obviously we have $(\Psi^{\#},\Phi^{\#})_n=I_c$, and
equation (\ref{cmeqsparam}) transforms into
\begin{equation}
\dot{x}=B^{\#}x+\Psi^{\#}(0)[{\cal L}(\alpha)-{\cal L}_0](\Phi^{\#}+
h(\alpha)U)x,
\label{cmeqsparamtrans}
\end{equation}
where $B^{\#}=U^{-1}BU$.  Let $T_U$ be the linear invertible
transformation of $\mbox{\rm Mat}_{c\times c}$ defined by $T_U(M)=U^{-1}MU$, 
let ${\cal T}^{\#}:\mbox{\rm Mat}_{c\times c}\longrightarrow\mbox{\rm Mat}_{c\times c}$
be defined by ${\cal T}^{\#}(M)=[B^{\#},M]$, and let ${\cal R}(\Psi^{\#}(0))$ be
the set of all $c\times c$ matrices whose columns are in the range of the matrix
$\Psi^{\#}(0)$.  Then it is easy to show that
\[
\mbox{\rm range}({\cal T}^{\#})=T_U(\mbox{\rm range}({\cal T}))
\]
and
\[
{\cal R}(\Psi^{\#}(0))=T_U({\cal R}(\Psi(0))).
\]
Consequently, from (\ref{homoleq}), we have that
\[
\mbox{\rm Mat}_{c\times c}=\mbox{\rm range}({\cal T}^{\#})\oplus \widehat{\cal W}^{\#},
\]
where $\widehat{\cal W}^{\#}=T_U(\widehat{\cal W})$.
Let $\{\Omega_1,\ldots,\Omega_{\delta}\}$ be a basis for ${\cal W}$ (see
(\ref{splitmatcc})), and ${\cal E}$ as in
Theorem~\ref{versunfconstthm}.  Let $A^m_0,\ldots,A^m_{c-q}$ be $n\times n$ matrices
such that
\begin{equation}
{\cal E}(\Omega_m)=\sum_{j=0}^{c-q}\,\Psi(0)A^m_j\Phi(\tau_j),\,\,\,\,\,m=1,\ldots,\delta.
\label{specbase}
\end{equation}
Then
\[
T_U({\cal E}(\Omega_m))=U^{-1}{\cal E}(\Omega_m)U=
\sum_{j=0}^{c-q}\,U^{-1}\Psi(0)A^m_j\Phi(\tau_j)U=
\sum_{j=0}^{c-q}\,\Psi^{\#}(0)A^m_j\Phi^{\#}(\tau_j),\,\,\,m=1,\ldots,\delta,
\]
i.e. given the delay times $\tau_0,\ldots,\tau_{c-q}$ which are such that
Lemma \ref{lemFM} holds, the matrices $A^m_0,\ldots,A^m_{c-q}$ which solve 
(\ref{specbase}) are such that the parametrized family (\ref{Lalphaconst}) 
generates a $\Lambda$-mini-versal unfolding of (\ref{fde})
independently 
of the choice of bases matrices $\Phi$ and $\Psi$ for $P$ and $P^*$ 
respectively, provided that $(\Psi,\Phi)_n=I_c$.

\subsection{Decomplexification\label{sec:decomp}}

In applications, it is usually the case that ${\cal L}_0$ in (\ref{fde}) and
${\cal L}(\alpha)$ in (\ref{fdep}) are real, i.e. they are
bounded linear operators from $C([-\tau,0],{\mathbb R}^n)$ into
${\mathbb R}^n$.  Although the previous theory has been carried out in 
complex spaces, it is straightforward to construct real versal unfoldings 
by a simple process of decomplexification of (\ref{Lalphaconst}).
We need to assume that the set $\Lambda$ defined in Section 2 is invariant 
under complex conjugation (one important example where this is always the 
case is center manifold reduction, in which case $\Lambda$ is the set of 
all roots of (\ref{char_eq}) with zero real parts).

\begin{thm}
Suppose that $\Lambda=\{\Lambda_0,\Lambda_{h},\ov{\Lambda_{h}}\}$
where $\Lambda_0$ is a subset of real eigenvalues and $\Lambda_{h}$
a subset of nonreal eigenvalues. Then, a real $\Lambda$-mini-versal 
unfolding of~(\ref{fde}) is given by 
\begin{equation}
{\cal L}(\alpha)={\cal L}_{0}+\sum_{p=1}^{\delta_0} \alpha_p L_p
+\sum_{s=\delta_0+1}^{\delta_0+\delta_h}\left( 
\beta_{s} {\rm Re}(L_{s})+\beta_{s+\delta_{h}}
{\rm Im}(L_{s})\right)
\label{realunf}
\end{equation}
where $\alpha_{p}\in {\mathbb R}$ for $p=1,\ldots,\delta_0$,
$\beta_{s},\beta_{s+\delta_{h}}\in {\mathbb R}$ for 
$s=\delta_0+1,\ldots,\delta_0+\delta_h$, $L_p$ is a bounded linear
operator from $C([-\tau,0],{\mathbb R}^n)$ into
${\mathbb R}^n$ for $p=1,\ldots,\delta_0$, and $L_s$ is a bounded
linear operator from $C([-\tau,0],{\mathbb R}^n)$ into ${\mathbb
  C}^n$, 
for $s=\delta_0+1,\ldots,\delta_0+\delta_h$.
\label{lem:decomp}
\end{thm}

\proof The $c\times c$ matrix $B$ can be decomposed as 
$B=\mbox{\rm diag}(B^{0},B^{h},
\ov{B^{h}})$ with $c=c_0+2c_{h}$, where $B^{0}$ is the 
$c_0\times c_0$ diagonal block of 
real eigenvalues $\Lambda_0$ in real Jordan canonical form and $B^{h}$ 
is the $c_{h}\times c_{h}$ diagonal block of nonreal eigenvalues 
$\Lambda_{h}$ in complex Jordan canonical form.

We establish the following notation for the remainder of the proof.
We let $\tilde{\varphi}_1,\ldots,\tilde{\varphi}_{c_0}$ and 
$\varphi_1,\ldots,\varphi_{c_{h}}$ be bases for the generalized eigenspace
corresponding to the eigenvalues of $\Lambda_0$ and $\Lambda_{h}$ 
respectively, chosen so that the matrices $\Phi^{0}=(\tilde{\varphi}_1,
\ldots,\tilde{\varphi}_{c_0})$ and $\Phi^{h}=(\varphi_1,\ldots,
\varphi_{c_{h}})$ satisfy respectively 
$\dot\Phi^{0}=\Phi^{0}\,B^{0}$ and 
$\dot{\Phi}^h=\Phi^h\,B^h$. Consequently, if we set 
$\Phi=(\Phi^{0},\Phi^h,\overline{\Phi^h})$, then the columns 
of $\Phi$ form a basis for $P$, and we have $\dot{\Phi}=\Phi\,B$.
Now, let $\tilde{\psi}_{1}^{*},\ldots,\tilde{\psi}_{c_0}^{*}$ be real 
linearly independent functions in $P^*$ (corresponding to $\Lambda_0$)
and
$\psi^*_1,\ldots,\psi^*_{c_{h}}$ be non-real linearly independent
functions in $P^*$
(corresponding to $\Lambda_h$).
If we denote 
$\Psi^{0*}=\mbox{\rm col}(\tilde{\psi}^*_1,\ldots,\tilde{\psi}^*_{c_0})$ and 
$\Psi^{h*}=\mbox{\rm col}(\psi^*_1,\ldots,\psi^*_{c_h})$,
then $\Psi^*=\mbox{\rm col}(\Psi^{0*},\Psi^{h*},\overline{\Psi^{h*}})$ 
is a basis for $P^*$. 
Define $\Psi=(\Psi^*,\Phi)_n^{-1}\Psi^*$, then $(\Psi,\Phi)_n=I_c$.  
Moreover, a simple computation shows that 
$\Psi=\mbox{\rm col}(\Psi^{0},\Psi^h,\overline{\Psi^h})$,
where $\Psi^{0}$ is a $c_0\times n$ real matrix corresponding to
$\Lambda_0$ whose rows are linearly
independent, 
and $\Psi^h$ is a $c_{h}\times n$ non-real matrix corresponding to
$\Lambda_h$ whose rows are 
linearly independent.

We define three projections $\Pi_0,\Pi_1,\Pi_2$  as in Remark~\ref{rankpsilem}
where $\Pi_0$ corresponds to $\Psi^{0}$ and has all real components while 
$\Pi_1$ and $\Pi_2$ correspond to $\Psi^{h}$ and $\ov{\Psi}^h$.
Consequently, the mappings $\Pi_1$ and $\Pi_2$ are such that
\begin{equation}
\overline{\Pi_1}=\Pi_2.  
\label{conjpi}
\end{equation}

Now let $\{\Omega_1,\ldots,\Omega_{\delta}\}$ be the elements of the basis of
${\cal W}$ used in the proof of Lemma \ref{homolproplem2} and Proposition \ref{homolprop}.
Assume that these basis elements are ordered so that the set
$\{\Omega_1,\ldots,\Omega_{\delta_0}\}$ corresponds to the block $B^0$,
$\{\Omega_{\delta_0+1},\ldots,\Omega_{\delta_0+\delta_h}\}$ corresponds to
the block $B^h$, $\{\Omega_{\delta_0+1+\delta_h},\ldots,\Omega_{\delta_0+2\delta_h}\}$
corresponds to the block $\ov{B^h}$, with
{\small \[
\Omega_p=\left(\begin{array}{ccc}\Omega_p^0&0_{c_0\times c_h}&0_{c_0\times c_h}\\
0_{c_h\times c_0}&0_{c_h\times c_h}&0_{c_h\times c_h}\\
0_{c_h\times c_0}&0_{c_h\times c_h}&0_{c_h\times c_h}\end{array}\right),\;\;\;
p=1,\ldots,\delta_0,
\]
\[
\Omega_s=\left(\begin{array}{ccc}0_{c_0\times c_0}&0_{c_0\times c_h}&0_{c_0\times c_h}\\
0_{c_h\times c_0}&\Omega^h_s&0_{c_h\times c_h}\\
0_{c_h\times c_0}&0_{c_h\times c_h}&0_{c_h\times c_h}\end{array}\right),\;\;
\Omega_{s+\delta_h}=\left(\begin{array}{ccc}0_{c_0\times c_0}&0_{c_0\times c_h}&0_{c_0\times c_h}\\
0_{c_h\times c_0}&0_{c_h\times c_h}&0_{c_h\times c_h}\\
0_{c_h\times c_0}&0_{c_h\times c_h}&\Omega^h_s\end{array}
\right),\;s=\delta_0+1,\ldots,\delta_0+\delta_h,
\]}
where $0_{k\times\ell}$ is the $k\times\ell$ zero matrix,
$\Omega_p^0$ is a $c_0\times c_0$ matrix with only one non-zero element, and
$\Omega^h_s$ is a $c_h\times c_h$ matrix with only one non-zero element.
Define
$R_{p}=(R_{p}^{0},0,0)$ for $p=1,\ldots,\delta_0$ where $0$ is the
$n\times c_{h}$ zero matrix and $R_{p}^{0}$ is a $n\times c_0$ matrix
with only one non-zero column, corresponding to a vector $v_{l}\in {\mathbb R}^{n}$ 
chosen as in Proposition~\ref{homolprop}, such that $\Psi(0)R_{p}$ 
has a $1$ at the same position as the $1$ in $\Omega_{p}$. Note that 
$v_{l}$ can be chosen in ${\mathbb R}^{n}$ since $\Pi_0$ is real. 
Now, for all $p=1,\ldots,\delta_0$, by Proposition~\ref{unfcandprop} we 
can find matrices $A_1^{p},\ldots,A_{c-q}^{p}$ such that
\[
(R_{p}^{0},0,0)=\sum_{j=0}^{c-q} A_{j}^{p} \,(\Phi_0(\tau_j),
\Phi^{h}(\tau_j),\ov{\Phi^{h}(\tau_j)}),
\]
where we have written $\Phi(\tau_{j})
\equiv(\Phi^0(\tau_j),\Phi^h(\tau_j),\overline{\Phi^h(\tau_j)})$.
The matrices $A_{j}^{p}$ can be chosen to be real matrices for all 
$p=1,\ldots,\delta_0$.  We then let $L_p$ be the bounded
linear operator from $C([-\tau,0],{\mathbb R}^n)$ into
${\mathbb R}^n$ defined by
$L_p(z)=\sum_{j=0}^{c-q} A_j^p\,z(\tau_j)$ for $p=1,\ldots,\delta_0$.

Now, define $R_s=(0,R^h_s,0)$ and $R_{s+\delta_h}=(0,0,\ov{R^h_s})$
for $s=\delta_0+1,\ldots,\delta_0+\delta_h$,
where the first $0$ is the $n\times c_0$ zero matrix and the second $0$ 
designates the $n\times c_{h}$ zero matrix.  The
$n\times c_h$ matrix $R^h_s$ has only one non-zero column,
corresponding to a vector $v_{l}\in {\mathbb C}^n$ chosen as in
Proposition \ref{homolprop}, such that $\Psi(0)R_s$ has a one
at the same position as the 1 in $\Omega_s$.  It then follows
(from (\ref{conjpi})) that $\Psi(0)\ov{R_s}$ has a one at the same
position as the 1 in ${\Omega_{s+\delta_h}}$.
Then, we set 
${\cal E}(\Omega_s)=\Psi(0)R_s$ and ${\cal E}(\Omega_{s+\delta_h})
=\Psi(0)R_{s+\delta_h}$.  From
Proposition \ref{unfcandprop}, there exist matrices 
$A^s_0,\ldots,A^s_{c-q}$ such that
\[
(0,R^h_{s},0)=\sum_{j=0}^{c-q}\,A^s_j \,(\Phi^{0}(\tau_{j}),
\Phi^h(\tau_j),\overline{\Phi^h(\tau_j)}),
\]
which we rewrite as the system
\[\begin{array}{lll}
{\dps\sum_{j=0}^{c-q}\,A^s_j\,\Phi^0(\tau_j)}=0, &\quad
{\dps\sum_{j=0}^{c-q}\,A^s_j\,\Phi^h(\tau_j)}=R^h_{s}, &\quad
{\dps\sum_{j=0}^{c-q}\,A^s_j\,\overline{\Phi^h(\tau_j)}}=0.
\end{array}
\]
If we take complex conjugates of the above system, we get
\[\begin{array}{lll}
{\dps\sum_{j=0}^{c-q}\,\overline{A^s_j}\Phi^0(\tau_j)=0}, &\quad
{\dps\sum_{j=0}^{c-q}\,\overline{A^s_j}\,\overline{\Phi^h(\tau_j)}}
=\overline{R^h_{s}}, & \quad
{\dps\sum_{j=0}^{c-q}\,\overline{A^s_j}\,{\Phi^h(\tau_j)}}=0,
\end{array}
\]
which is equivalent to
\[
(0,0,\overline{R^h_{s}})=\sum_{j=0}^{c-q}\,\overline{A^s_j} \,(
\Phi^{0}(\tau_{j}),\Phi^h(\tau_j),\overline{\Phi^h(\tau_j)}).
\]
Thus,
\[
R_{s+\delta_h}=\sum_{j=0}^{c-q}\,\overline{A^s_j}\Phi(\tau_j).
\]
From this, we conclude that the matrices $A^{s+{\delta}_{h}}_0,
\ldots,A^{s+{\delta}_{h}}_{c-q}$ of the decomposition
\[
R_{s+{\delta}_{h}}=\sum_{j=0}^{c-q}\,A^{s+{\delta}_{h}}_j\Phi(\tau_j)
\]
can be chosen so that $A^{s+{\delta}_{h}}_j=\overline{A^s_j}$,
$s=\delta_0+1,\ldots,
\delta_0+{\delta}_{h}$, $j=0,\ldots,c-q$.
Now, for all
$s\in\{\delta_0+1,\ldots,\delta_0+{\delta}_{h}\}$, we 
let $L_s$ and $L_{s+\delta_h}$ be the bounded linear operators from
$C([-\tau,0],{\mathbb R}^n)$ into ${\mathbb C}^n$
defined by
$L_s(z)=\sum_{j=0}^{c-q}\,A_j^s\,z(\tau_j)$ and
$L_{s+\delta_h}(z)=\sum_{j=0}^{c-q}\,\ov{A_j^s}\,z(\tau_j)$.
Taking real and imaginary parts of $L_{s}$ as in~(\ref{realunf})
yields a real $\Lambda$-mini-versal unfolding.
\qed

\begin{cor}
If $\Lambda=\Lambda_0$, then a real $\Lambda$-versal unfolding of~(\ref{fde}) 
is given by~(\ref{realunf}) with $\beta_{s}=\beta_{s+\delta_{h}}=0$  
for all $s=\delta_0+1,\ldots,\delta_0+\delta_{h}$.
\label{cor:realunf}
\end{cor}

\Section{Examples}
In this section, we illustrate our theory with several examples.
\begin{examp}
{\rm
Consider the scalar delay differential equation
\begin{equation}
\dot{x}(t)={\cal L}_0(x_t)=x(t)-x(t-1),
\label{ex1eq}
\end{equation}
which has been studied extensively \cite{FMTB,RLL}. The characteristic
equation has a double zero eigenvalue, therefore the center eigenspace
$P$ is two dimensional, i.e. $c=2$. 
A basis for $P$ is given by $\Phi=(1,\theta)$, and
we
notice that
\[
\mbox{\rm rank}(\mbox{\rm col}(\Phi(0),\Phi(-1)))=
\mbox{\rm
  rank}\,\left(\begin{array}{cc}1&0\\1&-1\end{array}\right)=2=c.
\]
By~(\ref{dimunfold}), $\delta=\dim({\cal W})=2$, 
and if follows
by Theorem~\ref{thm:scalar} and Corollary~\ref{cor:realunf} that a real 
$\Lambda$-mini-versal unfolding of (\ref{ex1eq}) is given by
\[
\dot{x}(t)={\cal L}(\alpha)(x_t)=x(t)-x(t-1)+\alpha_1x(t)+\alpha_2x(t-1)
\]
where $\alpha_1,\alpha_2\in {\mathbb R}$.
}
\end{examp}

\begin{examp}
{\rm Consider now the first order scalar equation
\begin{equation}\label{sc-eq}
\dot x={\cal L}_{0}x_{t}=A_1 x(t-\tau_1)+A_2 x(t-\tau_2)
\end{equation}
where $A_1,A_2$ are real parameters and the delays $\tau_1,\tau_2$ are
positive. B\'elair and Campbell~\cite{BC94} have shown that~(\ref{sc-eq}) 
has points of nonresonant double Hopf bifurcation in parameter space 
$(A_1,A_2,\tau_1,\tau_2)$ with eigenvalues $\pm \omega_1 i$ and 
$\pm\omega_2 i$. The center eigenspace $P$ is four-dimensional, thus
$c=4$. 

Suppose that~(\ref{sc-eq}) has a double Hopf bifurcation at 
$(A_1^{*},A_2^{*},\tau_1^{*},\tau_2^{*})$, with $\tau_1^*\neq\tau_2^*$.
By equation~(\ref{dimunfold}) we have that
$\delta=4$, and Theorem~\ref{thm:scalar} 
implies that a complex $\Lambda$-mini-versal unfolding
of~(\ref{sc-eq}) has
the form
\[
\dot x={\cal L}(\alpha)x_{t}=(A_1^{*}+\alpha_1) x(t-\tau_1^{*})
+(A_2^{*}+\alpha_2) x(t-\tau_2^{*})+\alpha_{3} x(t-\tau_{3})
+\alpha_4 x(t-\tau_4)
\]
for suitable choice of $\tau_3$ and $\tau_4$,
where $\alpha_{j}\in {\mathbb C}$ for $j=1,2,3,4$. 
}
\end{examp}

\begin{examp}
{\rm The second order scalar delay equation
\[
\ddot{x}(t)+\alpha\,\dot{x}(t)+\beta\,x(t)=f(x(t-\tau)),
\]
where $f$ is a smooth function such that $f(0)=0$,
was studied in Longtin and Milton~\cite{LM} as a model for the pupil light 
reflex. In Campbell and LeBlanc~\cite{C-LB}, it was shown that the
linearization of this equation about the trivial equilibrium, which we 
write in first-order form as
\begin{equation}
\begin{array}{lll}
\dot{x}_1(t)&=&x_2(t)\\
\dot{x}_2(t)&=&-\beta\,x_1(t)-\alpha\,x_2(t)+A\,x_1(t-\tau)
\end{array}
\label{CLeqlin}
\end{equation}
has a double Hopf point with $1:2$ resonance 
(eigenvalues $\pm i$ and $\pm 2i$) at parameter value
\begin{equation}
(\alpha,\beta,\tau,A)=\left(\,0,\,\frac{5}{2},\,\pi,\,\frac{-3}{2}\,\right).
\label{12parval}
\end{equation}
Using Theorem \ref{thmsuff}, it is straightforward to show that
if we treat $\alpha$, $\beta$, $\tau$ and $A$ as unfolding
parameters which
vary in a neighborhood of the point (\ref{12parval}), then
(\ref{CLeqlin}) generates a complex $\Lambda$-mini-versal unfolding for the 
singularity at parameter value
(\ref{12parval}).  Proving this 
amounts to constructing the $20\times 16$ matrix $S$ in
Theorem \ref{thmsuff} and showing that its rank is 16.

However, the goal of this example is to 
compute the real $\Lambda$-mini-versal unfolding 
for the singularity (\ref{CLeqlin}) (at parameter value (\ref{12parval})),
which results from
using the decomplexification procedure described in Section \ref{sec:decomp}.

We have that a complex basis for the center subspace is given by
\[
\Phi(\theta)=\left(\begin{array}{cccc}
e^{i\theta}&e^{2i\theta}&e^{-i\theta}&e^{-2i\theta}\\
ie^{i\theta}&2ie^{2i\theta}&-ie^{-i\theta}&-2ie^{-2i\theta}\end{array}
\right).
\]
A basis for the adjoint problem is given by
\[
\Psi^*(s)=\left(\begin{array}{lr}
e^{-is}&ie^{-is}\\
e^{-2is}&2ie^{-2is}\\
e^{is}&-ie^{is}\\
e^{2is}&-2ie^{2is}\end{array}
\right),
\]
which we renormalize by defining
\[
\Psi=(\Psi^*,\Phi)^{-1}\Psi^*.
\]
The result is such that
\begin{equation}
\Psi(0)=\kappa\left(\begin{array}{lr}
-2\,i \left( 3\,\pi +4\,i \right)  \left( -3\,\pi +4\,i \right) ^{2}
\,\,\,\,\,\,\,&
-32\,\pi +128\,i+6\,{\pi }^{3}-8\,i{\pi }^{2}\\[0.1in]
-i \left( -3\,\pi +8\,i \right)  \left( 3\,\pi +8\,i \right) ^{2}&
64\,\pi -16\,i{\pi }^{2}+64\,i-6\,{\pi }^{3}\\[0.1in]
2\,i \left( 3\,\pi -4\,i \right)  \left( -3\,\pi -4\,i \right) ^{2}&
-32\,\pi -128\,i+6\,{\pi }^{3}+8\,i{\pi }^{2}\\[0.1in]
i \left( -3\,\pi -8\,i \right)  \left( 3\,\pi -8\,i \right) ^{2}&
64\,\pi +16\,i{\pi }^{2}-64\,i-6\,{\pi }^{3}
\end{array}
\right),
\label{psi0cl}
\end{equation}
where $\kappa=(9\pi^4-32\pi^2-256)^{-1}$.
The matrix $B$ is given by
\[
\left(
\begin{array}{rrrr}
i&0&0&0\\
0&2i&0&0\\
0&0&-i&0\\
0&0&0&-2i
\end{array}
\right).
\]
A basis for ${\cal W}$ is given by the following four matrices
\[
\Omega_1=\left(\begin{array}{cccc}1&0&0&0\\0&0&0&0\\0&0&0&0\\0&0&0&0
\end{array}\right),\,\,\,\,\,\,\,\,\,\,\,\,\,
\Omega_2=\left(\begin{array}{cccc}0&0&0&0\\0&1&0&0\\0&0&0&0\\0&0&0&0
\end{array}\right),
\]
\[
\Omega_3=\left(\begin{array}{cccc}0&0&0&0\\0&0&0&0\\0&0&1&0\\0&0&0&0
\end{array}\right),\,\,\,\,\,\,\,\,\,\,\,\,\,
\Omega_4=\left(\begin{array}{cccc}0&0&0&0\\0&0&0&0\\0&0&0&0\\0&0&0&1
\end{array}\right).
\]
For any $\tau_0$ and $\tau_1$ such that $\tau_0\neq\tau_1$, the
matrix $\mbox{\rm col}(\Phi(\tau_0),\Phi(\tau_1))$ has rank
4. Therefore we choose, for example, $\tau_0=0$ and $\tau_1=-\pi$. 
Following Remark \ref{rankpsilem}, we construct
$\Pi_j:{\mathbb C}^2\longrightarrow {\mathbb C}$, $j=1,\ldots,4$:
\[
\Pi_j(v)=\psi_{j1}(0)v_1+\psi_{j2}(0)v_2,\,\,\,\,j=1,\ldots,4
\]
where $\psi_{j\ell}(0)$ are the elements in the matrix
$\Psi(0)$ in (\ref{psi0cl}), and $v=(v_1,v_2)^T\in {\mathbb C}^2$.
Therefore, we define
\[
R_1=\left(\begin{array}{cccc}
1&0&0&0\\
\frac{1-\psi_{11}(0)}{\psi_{12}(0)}&0&0&0\end{array}\right),\,\,\,\,\,
R_2=\left(\begin{array}{cccc}
0&1&0&0\\
0&\frac{1-\psi_{21}(0)}{\psi_{22}(0)}&0&0\end{array}\right),
\]
\[
R_3=\left(\begin{array}{cccc}
0&0&1&0\\
0&0&\frac{1-\overline{\psi_{11}(0)}}{\overline{\psi_{12}(0)}}&0
\end{array}\right),\,\,\,\,\,
R_4=\left(\begin{array}{cccc}
0&0&0&1\\
0&0&0&\frac{1-\overline{\psi_{21}(0)}}{\overline{\psi_{22}(0)}}
\end{array}\right).
\]
A simple computation shows that
\[
R_1=\frac{1}{4}\left(\begin{array}{cc}
1&-i\\
\frac{1-\psi_{11}(0)}{\psi_{12}(0)}&
-i\,\frac{1-\psi_{11}(0)}{\psi_{12}(0)}
\end{array}\right)\Phi(0)+
\frac{1}{4}\left(\begin{array}{cc}
-1&i\\
-\frac{1-\psi_{11}(0)}{\psi_{12}(0)}&
i\,\frac{1-\psi_{11}(0)}{\psi_{12}(0)}\end{array}\right)\Phi(-\pi)
\]
and
\[
R_2=\frac{1}{8}\left(\begin{array}{cc}
2&-i\\
\frac{2(1-\psi_{21}(0))}{\psi_{22}(0)}&
-2i\,\frac{1-\psi_{21}(0)}{\psi_{22}(0)}
\end{array}\right)\Phi(0)+
\frac{1}{8}\left(\begin{array}{cc}
2&-i\\
\frac{2(1-\psi_{21}(0))}{\psi_{22}(0)}&
-2i\,\frac{1-\psi_{21}(0)}{\psi_{22}(0)}
\end{array}\right)\Phi(-\pi).
\]
Thus, if we set
\[
\beta_1+i\gamma_1=\frac{1-\psi_{11}(0)}{\psi_{12}(0)},\,\,\,\,\,\,
\beta_2+i\gamma_2=\frac{1-\psi_{21}(0)}{\psi_{22}(0)},
\]
then the operators $L_m$ in (\ref{Lmdef}) are
\[
L_1(z)=\frac{1}{4}\left(\begin{array}{cc}
1&-i\\
\beta_1+i\gamma_1&\gamma_1-i\beta_1\end{array}\right)z(0)+
\frac{1}{4}\left(\begin{array}{cc}
-1&i\\
-\beta_1-i\gamma_1&-\gamma_1+i\beta_1\end{array}\right)z(-\pi),
\]
\[
L_2(z)=\frac{1}{8}\left(\begin{array}{cc}
2&-i\\
2\beta_2+2i\gamma_2&2\gamma_2-2i\beta_2\end{array}\right)z(0)+
\frac{1}{8}\left(\begin{array}{cc}
2&-i\\
2\beta_2+2i\gamma_2&2\gamma_2-2i\beta_2\end{array}\right)z(-\pi),
\]
\[
L_3=\overline{L_1},\,\,\,\,\,L_4=\overline{L_2}.
\]
Thus, the decomplexification procedure yields the
following real $\Lambda$-mini-versal unfolding of the singularity 
(\ref{CLeqlin}) (at parameter value (\ref{12parval}))
\[
\begin{array}{lll}
\left(\begin{array}{c}
\dot{x}_1(t)\\\dot{x}_2(t)\end{array}\right)&=&
\left(\begin{array}{c}
x_2(t)\\
-\frac{5}{2}\,x_1(t)-\frac{3}{2}\,x_1(t-\pi)
\end{array}\right)
+\\[0.2in]
&&\left(\alpha_1\left(\begin{array}{cc}1&0\\\beta_1&\gamma_1\end{array}\right)
+\alpha_2\left(\begin{array}{cc}0&-1\\\gamma_1&-\beta_1\end{array}\right)\right)
\left(\begin{array}{c}
x_1(t)-x_1(t-\pi)\\
x_2(t)-x_2(t-\pi)\end{array}\right)+\\[0.2in]
&&\left(\alpha_3\left(\begin{array}{cc}2&0\\2\beta_2&2\gamma_2\end{array}\right)+\alpha_4\left(\begin{array}{cc}0&-1\\2\gamma_2&-2\beta_2\end{array}\right)\right)
\left(\begin{array}{c}
x_1(t)+x_1(t-\pi)\\
x_2(t)+x_2(t-\pi)\end{array}\right).
\end{array}
\]
}
\end{examp}

\vspace*{0.2in}
\noindent
{\Large\bf Appendix}

\appendix
\Section{Proof of Lemma~\ref{rangeTlem}}
We use the following inner product on the space 
$\mbox{\rm Mat}_{c\times c}$:
\begin{equation}
\langle\,M_1,M_2\,\rangle = \mbox{\rm trace}(M_1M_2^*).
\label{inprod}
\end{equation}
A simple computation shows that with respect to
(\ref{inprod}), the conjugate transpose of ${\cal T}$ is given by
${\cal T}^*(M)=[B^*,M]$.  From the Fredholm alternative, we have
$\mbox{\rm range}({\cal T})=\mbox{\rm ker}({\cal T}^*)^\perp$.
Recall that $\mbox{\rm ker}({\cal T}^*)$ has been characterized in Lemma
\ref{commstruct}.
So $Y\in\mbox{\rm ker}({\cal T}^*)^\perp$, if and
only if $Y$ is orthogonal to
all elements in a basis of $\mbox{\rm ker}({\cal T}^*)$.
We construct a basis for $\mbox{\rm ker}({\cal T}^*)$ as follows:
for $j\in\{1,\ldots,r\}$, denote
\[
n_o(j)=\mbox{\rm number of distinct oblique segments in the
block ${\cal M}_j$, as in Figure \ref{fig1}},
\]
and give some ordering to these oblique segments, numbering them from 1 to
$n_o(j)$.
Then, for $j\in\{1,\ldots,r\}$ and $\ell\in\{1,\ldots,n_o(j)\}$, we
define
$M_{j,\ell}$ to be the matrix structured as in (\ref{Ndiagstruct}),
and such that all diagonal blocks are zero except the block $j$;
and in this block $j$, the only oblique segment which is non-zero
is the $\ell^{th}$ segment, whose elements all have
value 1.  Then $\{M_{j,\ell}\}$ forms a basis for
$\mbox{\rm ker}({\cal T}^*)$.

Now consider $Y\in\mbox{\rm Mat}_{c\times c}$, which we
partition as in (\ref{rangegen}).  Multiply $Y$ on the right
by $M_{j,\ell}^*$, for
some $j\in\{1,\ldots,r\}$, $\ell\in\{1,\ldots,n_o(j)\}$.
Then $YM_{j,\ell}^*$ is partitioned as in 
(\ref{rangegen}), and is such that
\[
\begin{array}{ccc}
YM_{j,\ell}^*&=&\left(
\begin{array}{ccccccc}
0&\cdots&0&{\cal Y}_{1,j}{\cal M}_{j,\ell}^*&0&\cdots&0\\
0&\cdots&0&{\cal Y}_{2,j}{\cal M}_{j,\ell}^*&0&\cdots&0\\
\vdots&\vdots&\vdots&\vdots&\vdots&\vdots&\vdots\\
\vdots&\vdots&\vdots&\vdots&\vdots&\vdots&\vdots\\
\vdots&\vdots&\vdots&\vdots&\vdots&\vdots&\vdots\\
\vdots&\vdots&\vdots&\vdots&\vdots&\vdots&\vdots\\
0&\cdots&0&{\cal Y}_{r,j}{\cal M}_{j,\ell}^*&0&\cdots&0
\end{array}
\right),\\
&&\uparrow\\
&&\mbox{\rm $j^{th}$ vertical block}\\
&&\mbox{\rm in partitioning (\ref{rangegen})}
\end{array}
\]
where ${\cal M}_{j,\ell}$ is the $j^{th}$ diagonal block
in the decomposition (\ref{Ndiagstruct}) of $M_{j,\ell}$.
Thus, $\langle\,Y,M_{j,\ell}\,\rangle=\mbox{\rm trace}(YM_{j,\ell}^*)=
\mbox{\rm trace}({\cal Y}_{j,j}{\cal M}_{j,\ell}^*)$.
It follows that $Y\in\mbox{\rm ker}({\cal T}^*)^\perp$ if and only if 
for all $j\in\{1,\ldots,r\}$, the
block ${\cal Y}_{j,j}$ in the partition (\ref{rangegen}) is such that
$\mbox{\rm trace}({\cal Y}_{j,j}{\cal M}_{j,\ell}^*)=0$
for all $\ell\in\{1,\ldots,n_o(j)\}$.
Using the partition illustrated in Figure \ref{fig1}, we write
\begin{equation}
{\cal Y}_{j,j}=\left(\begin{array}{ccccc}
{\cal Y}_{j,j}^{1,1}&{\cal Y}_{j,j}^{1,2}&\cdots&\cdots&{\cal
  Y}_{j,j}^{1,k_j}\\
{\cal Y}_{j,j}^{2,1}&{\cal Y}_{j,j}^{2,2}&\cdots&\cdots&{\cal
  Y}_{j,j}^{2,k_j}\\
\vdots&\vdots&\vdots&\vdots&\vdots\\
\vdots&\vdots&\vdots&\vdots&\vdots\\
{\cal Y}_{j,j}^{k_j,1}&{\cal Y}_{j,j}^{k_j,2}&\cdots&\cdots&{\cal
  Y}_{j,j}^{k_j,k_j}
\end{array}
\right)
\label{Adecomp}
\end{equation}
and
\begin{equation}
{\cal M}_{j,\ell}=\left(\begin{array}{ccccc}
{\cal M}_{j,\ell}^{1,1}&{\cal M}_{j,\ell}^{1,2}&\cdots&\cdots&{\cal
  M}_{j,\ell}^{1,k_j}\\
{\cal M}_{j,\ell}^{2,1}&{\cal M}_{j,\ell}^{2,2}&\cdots&\cdots&{\cal
  M}_{j,\ell}^{2,k_j}\\
\vdots&\vdots&\vdots&\vdots&\vdots\\
\vdots&\vdots&\vdots&\vdots&\vdots\\
{\cal M}_{j,\ell}^{k_j,1}&{\cal M}_{j,\ell}^{k_j,2}&\cdots&\cdots&{\cal
  M}_{j,\ell}^{k_j,k_j}
\end{array}
\right),
\label{basdecomp}
\end{equation}
where ${\cal Y}_{j,j}^{p,q}$ and ${\cal M}_{j,\ell}^{p,q}$ are
$n_{j,p}\times n_{j,q}$, for $p,q\in\{1,\ldots,k_j\}$.
Let $\tilde{p}$ and $\tilde{q}$ be such that 
${\cal M}_{j,\ell}^{\tilde{p},\tilde{q}}$ is
the unique block in (\ref{basdecomp}) which has the non-zero oblique segment.
Then ${\cal Y}_{j,j}{\cal M}_{j,\ell}^*$ is partitioned as in
(\ref{Adecomp}) and has the form
\[
\begin{array}{ccc}
{\cal Y}_{j,j}{\cal M}_{j,\ell}^*&=&\left(\begin{array}{ccccccc}
0&\cdots&0&{\cal Y}_{j,j}^{1,\tilde{q}}\cdot({\cal
  M}_{j,\ell}^{\tilde{p},\tilde{q}})^*&0&\cdots&0\\
0&\cdots&0&{\cal Y}_{j,j}^{2,\tilde{q}}\cdot({\cal
  M}_{j,\ell}^{\tilde{p},\tilde{q}})^*&0&\cdots&0\\
\vdots&\vdots&\vdots&\vdots&\vdots&\vdots&\vdots\\
\vdots&\vdots&\vdots&\vdots&\vdots&\vdots&\vdots\\
0&\cdots&0&{\cal Y}_{j,j}^{k_j,\tilde{q}}\cdot({\cal
  M}_{j,\ell}^{\tilde{p},\tilde{q}})^*&0&\cdots&0
\end{array}
\right).\\
&&\uparrow\\
&&\mbox{\rm $\tilde{p}^{th}$ vertical block in}\\
&&\mbox{\rm partitioning (\ref{Adecomp})}
\end{array}
\]
It follows that 
$Y\in\mbox{\rm ker}({\cal T}^*)^\perp$ if and only if
$Y$ is of the form (\ref{rangegen}), and
for all $j\in\{1,\ldots,r\}$ and for all $\ell\in\{1,\ldots,n_o(j)\}$,
we have
$\mbox{\rm trace}({\cal Y}_{j,j}^{\tilde{p},\tilde{q}}\cdot
({\cal M}_{j,\ell}^{\tilde{p},\tilde{q}})^*)=0$, where 
$\tilde{p}$ and $\tilde{q}$ are such that 
${\cal M}_{j,\ell}^{\tilde{p},\tilde{q}}$ is
the unique block in (\ref{basdecomp}) which has the non-zero oblique
segment.
Now, suppose the block
${\cal M}_{j,\ell}^{\tilde{p},\tilde{q}}$ is such that
the oblique segment whose elements are all 1's is the one
which passes through the $\nu^{th}$ leftmost element
in the bottom row of the block, where
$\nu\in\{1,\ldots,n_{j,\tilde{q}}\}$.
Then a simple computation shows that
${\cal Y}_{j,j}^{\tilde{p},\tilde{q}}\cdot
({\cal M}_{j,\ell}^{\tilde{p},\tilde{q}})^*$ is an 
$n_{j,\tilde{p}}\times n_{j,\tilde{p}}$ matrix
whose first $n_{j,\tilde{p}}-\nu$ rows are zero, and the remaining rows are
the rows $1,\ldots,\nu$ of ${\cal Y}_{j,j}^{\tilde{p},\tilde{q}}$.
It follows that
$\mbox{\rm trace}({\cal Y}_{j,j}^{\tilde{p},\tilde{q}}\cdot
({\cal M}_{j,\ell}^{\tilde{p},\tilde{q}})^*)$ is the
sum of the elements in the oblique segment of the
block ${\cal Y}_{j,j}^{\tilde{p},\tilde{q}}$ which is the
same oblique segment as the non-zero segment of
${\cal M}_{j,\ell}^{\tilde{p},\tilde{q}}$, from which we get the
conclusion.
\hfill
\qed

\vspace*{0.25in}
\noindent
{\Large\bf Acknowledgements}

\vspace*{0.2in}
This research is partly supported by the
Natural Sciences and Engineering Research Council of Canada in the
form of a postdoctoral fellowship (PLB) and an individual research
grant (VGL).

\end{document}